%% file: main.tex
\documentclass[11pt,letterpaper]{article}
\input{preamble}

\input{macros}

\title{The limit shape of uniform random antichains in $[0,1]^d$}
\author{Andrew Elvey Price\thanks{CNRS, Université de Tours, Université d'Orléans, IDP, UMR 7013, F-37200 Tours, France.}\and Minjae Park\thanks{Department of Mathematics and Statistics, Auburn University, USA.}}
\date{\today}
\begin{document}
\maketitle

\begin{abstract}
      For each $d\ge2$, we study $n$ independent uniform points in $[0,1]^d$, conditioned to form an antichain: no point is strictly smaller than another in every coordinate. We prove that the point cloud converges in Hausdorff distance to the closure of the hypersurface
      \[
            \mathcal S_d=\left\{\mathbf x\in(0,1)^d: \sum_{i=1}^d\Phi^{-1}(\mathbf x_i)=0\right\}
      \]
      in probability as $n\to\infty$, where $\Phi$ is the standard normal cumulative distribution function. The proof also shows that the antichain probability $a_n$ satisfies $\log a_n=-\frac{n}{d-1}\log n+\beta_d n+o(n)$ for a finite constant $\beta_d$.
      We identify the limiting empirical distribution and hence the limiting permuton of uniform random $d$-dimensional permutations avoiding a single pattern of size $2$.
\end{abstract}

\begin{figure}[!htbp]
      \centering
      \includegraphics[width=0.47\linewidth]{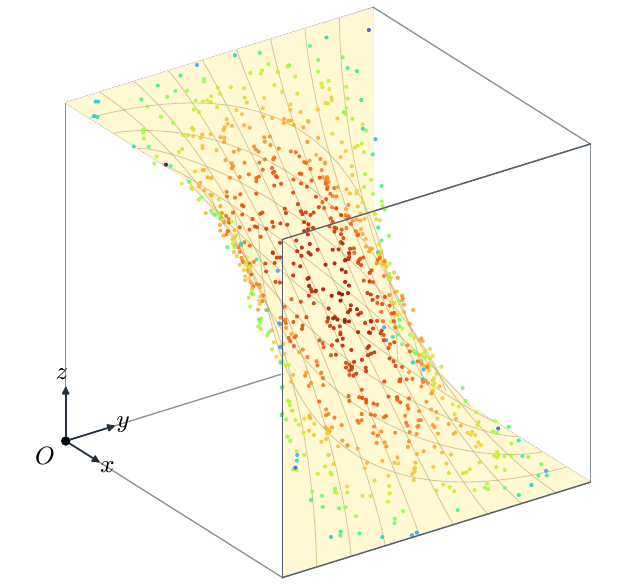}\qquad
      \includegraphics[width=0.47\linewidth]{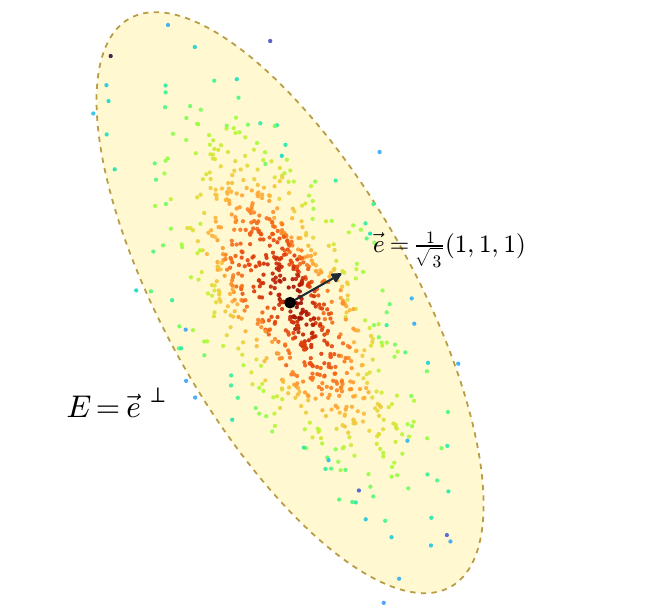}
      \caption{\textbf{Left:} The surface $\mathcal S_3=\{\mathbf x\in(0,1)^3: \Phi^{-1}(\mathbf x_1)+\Phi^{-1}(\mathbf x_2) +\Phi^{-1}(\mathbf x_3)=0\}$ and a uniform $900$-point antichain sampled by Markov chain Monte Carlo. \textbf{Right:} The same points in Gaussian coordinates $\mathbf g=\Phi^{-1}(\mathbf x)$, with $\Phi^{-1}$ applied coordinatewise, together with a truncated portion of the plane $E=\{\mathbf g\in\R^3: \mathbf g_1+\mathbf g_2+\mathbf g_3=0\}$. Colors indicate $\|\mathbf g\|$ in both coordinate systems.}
      \label{fig:limit-shape}
\end{figure}

\newpage

\tableofcontents

\section{Introduction}

For two points $\vx=(\vx_{1},\ldots,\vx_{d})\in\R^d$ and $\vy=(\vy_{1},\ldots,\vy_{d})\in\R^d$, we write $\vx\prec\vy$ if $\vx_{i}<\vy_{i}$ for each $i\in\{1,\ldots,d\}$. We call a set of points in $\R^d$ an \emph{antichain} if it contains no distinct points $\vx,\vy$ with $\vx\prec\vy$, that is, it is an antichain with respect to the order $\prec$. We study the behavior of uniform random antichains of $n$ points in the cube $[0,1]^d$ as $n\to\infty$. Our main results determine the limit shape of a random antichain and the distribution over this shape. We also improve on previous results regarding the asymptotic form of the probability that $n$ uniform random points form an antichain.

Throughout this paper we fix the dimension $d\ge2$. Let $\vX^{(1)},\ldots,\vX^{(n)}$ be independent uniform points in $[0,1]^d$, and define $\cA_n$ to be the event that these points form an antichain, and $a_{n}$ to be the probability that this occurs:
\[
      \cA_n=\{\vX^{(k)}\not\prec\vX^{(\ell)}\text{ for all }k\ne\ell\}, \qquad a_n=\P(\cA_n).
\]
The first result establishes a linear-order correction to the leading asymptotic of $\log a_n$:
\begin{theorem}\label{thm:counting}
      There is a finite constant $\beta_d$ such that

      \begin{equation}\label{eq:counting}
            \log a_n=-\frac n{d-1}\log n+\beta_d n+o(n).
      \end{equation}
\end{theorem}
As we will discuss, the leading term of this asymptotic was already known. However, the fact that the remaining term, divided by $n$, converges to a finite constant is new for $d\ge3$. The result follows directly from Proposition~\ref{prop:volume}; see Remark~\ref{rem:proof-thm-1}.

Our other results involve the cumulative distribution function $\Phi$ of the standard normal distribution, and its density $\phi$, defined as follows:
\[
      \phi(s)=\frac1{\sqrt{2\pi}}e^{-\frac{s^2}{2}},\qquad \Phi(s)=\int_{-\infty}^s\phi(u)\dd u.
\]
In particular, we show in Theorem~\ref{thm:shape} that uniform random antichains typically approximate the hypersurface
\[
      \cS_d=\left\{\vx\in(0,1)^d:\sum_{i=1}^d\Phi^{-1}(\vx_i)=0\right\}.
\]
Writing $\ve=\frac{1}{\sqrt d}(1,\ldots,1)$ and $E=\ve^\perp$, the limit shape $\cS_d$ can be defined as $\cS_d=\Phi(E)$, where $\Phi$ acts coordinatewise.

To make this statement precise, we define the Hausdorff distance between nonempty compact subsets $A,B\subset[0,1]^d$ by
\[
      d_{\rm H}(A,B)=\max\left\{\sup_{\vx\in A}\dist(\vx,B), \sup_{\vy\in B}\dist(\vy,A)\right\}, \qquad \dist(\vx,B)=\inf_{\vy\in B}\|\vx-\vy\|,
\]
where closures and distances are taken in the closed Euclidean cube $[0,1]^d$.
The Hausdorff distance controls both outliers and holes. Let $A_n=\{\vX^{(1)},\ldots,\vX^{(n)}\}$ denote the uniformly sampled point set and write $\P_n=\P(\,\cdot\mid\cA_n)$ for probability under the uniform antichain model.
The following theorem identifies the shape filled by this set as $\overline{\cS_d}$.
See Figure~\ref{fig:limit-shape} for an illustration when $d=3$.

\begin{theorem}\label{thm:shape}
      For every $\delta>0$, there is some $\alpha=\alpha(d,\delta)>0$ such that
      \[
            \P_n\left(d_{\rm H}(A_n,\overline{\cS_d})>\delta \right) \le e^{-\alpha n}
      \]
      for all sufficiently large $n$. In particular, the entire sample fills $\overline{\cS_d}$.
\end{theorem}
The proof is given at the end of Section~\ref{sec:shortcuts}, after we identify the limiting distribution of points along this shape and establish empirical concentration. This concentration, together with the antichain condition, yields Theorem~\ref{thm:shape}.

\subsection{Related results and applications}\label{sec:background}
A $d$-permutation of $[n]=\{1,\ldots,n\}$ is a $d$-tuple of permutations $(\mathrm{id},\pi_{2},\ldots,\pi_{d})$ of $[n]$ with $\mathrm{id}$ denoting the identity permutation on $[n]$. There is a natural correspondence between samples $\vX^{(1)},\ldots,\vX^{(n)}\in[0,1]^{d}$ and $d$-permutations, since coordinate ties occur with probability zero. Sorting the points $\vX^{(1)},\ldots,\vX^{(n)}$ by their first coordinate and recording the remaining coordinate ranks gives a $d$-permutation $\pi=(\mathrm{id},\pi_2,\ldots,\pi_d)$ of $[n]$. In the case that the original sample is uniform, $\pi$ is a uniform $d$-permutation. Moreover, the original sample is an antichain if and only if there is no $k<\ell$ satisfying $\pi_i(k)<\pi_i(\ell)$ for every $i\ge2$, that is, if $\pi$ lies in the class $\cC_n$ of $d$-permutations avoiding the pattern $(12,12,\ldots,12)$, which we define below. Theorem~\ref{thm:counting} yields information about the number of these permutations, as $|\cC_n|=(n!)^{d-1}a_n$ (see Proposition~\ref{cor:counting}). A \emph{$d$-dimensional permuton} is a probability measure on $[0,1]^d$ with uniform coordinate marginals; such measures describe permutation limits~\cite{HKMRS,BL}. Theorem~\ref{thm:permuton} identifies the limiting permuton of $\cC_n$.

A $d$-permutation $\pi=(\mathrm{id},\pi_{2},\ldots,\pi_{d})$ of $[n]$ is said to contain a $d$-permutation $\tau=(\mathrm{id},\tau_{2},\ldots,\tau_{d})$ of $[m]$, with $m\le n$, as a pattern if there are indices $1\le i_{1}<\cdots<i_{m}\le n$ such that $(\pi_{j}(i_k))_{k\in[m]}$ has the same relative order as $(\tau_{j}(k))_{k\in[m]}$ for $j=2,\ldots,d$. Otherwise $\pi$ {\em avoids} $\tau$.
By symmetry, Theorems~\ref{thm:shape} and~\ref{thm:empirical-limit} determine the geometric and permuton limits of $d$-permutations avoiding a pattern of length $2$ (see Theorem~\ref{thm:permuton} and Remark~\ref{rem:rank-hausdorff}), while Theorem~\ref{thm:counting} determines the asymptotic enumeration of these $d$-permutations (see Proposition~\ref{cor:counting}).

The modern study of classical pattern-avoiding permutations, corresponding to $d=2$, grew from Knuth's work on permutation sorting in {\em The Art of Computer Programming}~\cite{Knuth} and has since been studied extensively; see~\cite{Vatter} for a survey. In particular, the Marcus--Tardos theorem~\cite{marcus2004excluded} shows that for any single pattern $\tau$, the number of $\tau$-avoiding permutations of length $n$ grows at most exponentially in $n$. In recent years, the study of limit permutons has attracted a lot of attention~\cite{borga2020local,garbe2024pattern,hohmeier2025permuton,bassino2018brownian}. Early examples of pattern-avoiding $3$-permutations arose in the consideration of transformations via priority queues~\cite{atkinson1995priority}. While a number of articles have addressed the enumeration of $d$-permutations avoiding multiple patterns~\cite{aldred2005permuting,EP2026stretched,BM}, single-pattern classes remain much less understood.

The Marcus--Tardos theorem was generalized to $d>2$ independently in~\cite[Theorem~1.9]{CK} and~\cite[Theorem~3.4]{GP}, which show that for any fixed single pattern $\tau$ of size at least $2$, the number $\text{Av}_{n}(\tau)$ of $d$-permutations of $[n]$ avoiding $\tau$ satisfies, for some $C=C(d,\tau)>0$,
\[
      e^{-C n}\le\text{Av}_{n}(\tau)\Big/(n!)^{d-1-\frac1{d-1}}\le e^{C n}.
\]
Question~9.3 of~\cite{CK} and the last section of~\cite{GP} ask whether the $n$th root of this ratio converges to a positive constant. In Proposition~\ref{cor:counting}, we answer this question in the affirmative in the case that $\tau$ has size $2$.

Permuton limits of $d$-permutations avoiding multiple patterns for $d>2$ have also been studied in~\cite{BL}. To our knowledge, for $d>2$, Theorem~\ref{thm:permuton} is the first determination of a permuton limit for a class avoiding a single pattern. We believe that this opens the door to further results; in particular, we expect that many different classes of $d$-permutations avoiding a single pattern share this same limit shape (see Conjecture~\ref{conj:bounded-chains} and Question~\ref{ques:universality}). The methods may also allow other $d$-permutation classes to be analyzed, whose permuton limits have different shapes. The Gaussian coordinate change $\Phi^{-1}$ offers a new geometric approach to such problems, and may even be useful for classes of permutations ($d=2$) that are so far not well understood.

Our model is also related to random multidimensional orders, which were introduced in~\cite{winkler1985random}. The random partial order $P_k(n)$ is defined as the intersection of $k$ independent uniform
linear orders on $[n]$, meaning that $a<b$ in this order if and only if $a<b$ in each of the $k$ linear orders. Let $e(P)$ count the linear
extensions of $P$, that is, the total linear orders respecting all
comparisons in $P$. Given uniform points $\vX^{(1)},\ldots,\vX^{(n)}\in [0,1]^{d}$, a partial order with the law of $P_{d-1}(n)$ can be defined by taking $a<b$ if and only if $\vX^{(a)}_{i}<\vX^{(b)}_{i}$ for each $i\in\{1,\ldots,d-1\}$. Since the reversed last coordinate order
extends this order if and only if $\vX^{(1)},\ldots,\vX^{(n)}$ form an antichain, we have
$a_n=\E[e(P_{d-1}(n))]/n!$. Indeed, applying the estimate~\cite[Corollary~4]{Brightwell} to the right-hand side yields the lower bound for $a_n$ obtained in~\cite[Theorem~8.4]{CK}.

In dimension three, $\cC_n$ is in bijection with intervals in the weak Bruhat order through the correspondence with avoidance of pairs of permutations~\cite{GP,BM}. This links the model to the comparability problems studied in~\cite{HP,BCHHLM}. Proposition~\ref{prop:intervals} identifies the joint limiting geometry of the two bounding permutations of an interval.
Their limiting inversion densities are also identified in Corollary~\ref{cor:interval-inversions}. Section~\ref{sec:applications} contains these applications and open questions. Related stepped surfaces have also been studied through surface tension~\cite{Sheffield}. For our model in dimension three, the surface tension was identified in~\cite{KenyonComm}, and our theorem confirms that the limit shape is the minimizer of the corresponding surface tension functional.

\subsection{Setting and proof strategy}

We start by discussing the relevance of the coordinatewise Gaussian transform $\Phi^{-1}:(0,1)^d\to\R^{d}$. The use of this transformation was suggested by OpenAI's ChatGPT (see~\nameref{sec:ai}), and we do not know how this idea arose exactly, but we give some discussion here of where it could have come from. Fundamentally, this transformation is useful because it separates a point cloud's location from its arrangement. Starting with $\mathbf x\in(0,1)^d$, apply a smooth bijection $Q:(0,1)\to\R$ with $Q'>0$ to each coordinate, setting $s_i=Q(x_i)$. This gives $\mathbf s\in\R^d$ and preserves the antichain condition. The change of variables gives $\dd\mathbf x=(\prod_iq(s_i))\,\dd\mathbf s$, where $q=(Q^{-1})'$. We can view this density as a Boltzmann weight $e^{-\sum_iV(s_i)}$, with $V=-\log q$ assigning an energy to each coordinate. The logarithm of the unweighted volume of available arrangements plays the role of counting entropy.

For $n$ transformed points $\mathbf s_1,\ldots,\mathbf s_n\in\R^d$, let $\mathbf m$ be their average and write $\mathbf s_j=\mathbf m+\mathbf u_j$. The centered points $\mathbf u_1,\ldots,\mathbf u_n$ retain all pairwise differences, and hence the antichain condition. Since $\sum_j\mathbf u_j=0$, the quadratic part of the Gaussian energy splits as
\[
      \frac12\sum_{j=1}^n\|\mathbf s_j\|^2=\frac n2\|\mathbf m\|^2+\frac12\sum_{j=1}^n\|\mathbf u_j\|^2.
\]
The first term depends only on the group's location, and the second only on its arrangement. Moving the barycenter therefore changes every centered arrangement's weight by the same factor, allowing us to study arrangements separately from the cost of their location. Among smooth positive densities on $\R$, requiring this separation for every configuration singles out the Gaussian family. Choosing an origin and scale then gives $Q=\Phi^{-1}$.

Accordingly, we set
\begin{equation}\label{eq:transform_Phi}
      \vG^{(k)}=\Phi^{-1}(\vX^{(k)})\qquad \text{for $k=1,\ldots,n$}.
\end{equation}
This is well defined almost surely and gives independent standard Gaussian vectors before conditioning, with $\cA_n$ unchanged. In these coordinates, the candidate limit shape $\cS_d=\Phi(E)$ becomes the hyperplane $E$. We now split each point $\vG^{(k)}$ into its {\em floor position} $\vF^{(k)}\in E$ and its (signed) {\em height} $H^{(k)}=\langle\vG^{(k)},\ve\rangle$:
\begin{equation}\label{eq:split}
      \vG^{(k)}=\vF^{(k)}+H^{(k)}\ve,\qquad \vF^{(k)}=\Proj_E\vG^{(k)}.
\end{equation}
Here $\Proj_E\vg=\vg-\langle\vg,\ve\rangle\ve$ is orthogonal projection onto $E$.
Since the Gaussian distribution on $\R^{d}$ is rotationally invariant, these floor positions and heights are independent standard Gaussians on $E$ and $\R$, respectively (before conditioning). On the other hand, after conditioning on the antichain event $\cA_{n}$, our results will imply that in the limit $n\to\infty$, the typical height $H^{(k)}$ converges to $0$, while the vectors $\vF^{(k)}$ follow a rescaled Gaussian distribution on $E$.

Our main Theorems~\ref{thm:counting}, \ref{thm:shape} and~\ref{thm:empirical-limit} are all proven using Proposition~\ref{prop:volume}, an asymptotic estimate for a weighted Gaussian expectation. For a bounded continuous floor weight $\sigma:E\to\R$ and a parameter $c\in\R$, define the (possibly infinite) weighted expectation
\begin{equation}\label{eq:def-zn}
      \cZ_n(\sigma,c) =\E\left[\ind_{\cA_n} \exp\biggl\{\sum_{k=1}^n\Bigl(\sigma(\vF^{(k)})+c(H^{(k)})^2\Bigr)\biggr\}\right].
\end{equation}
Here $\ind_{\cA_n}$ is $1$ on antichains and $0$ otherwise, so this expectation is an integral over labeled antichains in the cube. The factor $e^\sigma$ weights floor positions, while positive $c$ favors larger absolute heights. Boundedness of $\sigma$ makes the integral finite for $c<1/2$. In particular, $\cZ_n(0,0)=a_n$.

The density of a standard Gaussian floor position is
\begin{equation}\label{eq:rho}
      \psi_E(\vu) =(2\pi)^{-\frac{d-1}{2}}e^{-\frac{\|\vu\|^2}{2}},\qquad \vu\in E.
\end{equation}
Integration on $E$ uses its induced Euclidean measure, equivalently Lebesgue measure in any orthonormal coordinates on $E$. The following proposition gives the weighted estimate underlying Theorem~\ref{thm:counting}.

\begin{proposition}\label{prop:volume}\label{thm:volume}
      There is a finite constant $\gamma_d$ such that, for every bounded continuous $\sigma:E\to\R$ and every $c<1/2$,
      \begin{equation}\label{eq:weighted}
            \lim_{n\to\infty}\left(\frac1n\log\cZ_n(\sigma,c)+\frac{\log n}{d-1}\right) =\gamma_d-\frac12\log(2\pi)+\frac d{d-1}\log\int_E(\psi_E(\vu)e^{\sigma(\vu)})^{\frac{d-1}d}\dd\vu.
      \end{equation}
\end{proposition}

\begin{remark}\label{rem:proof-thm-1}
      The constant $\gamma_d$ is defined in~\eqref{eq:gamma}. Setting $(\sigma,c)=(0,0)$ yields Theorem~\ref{thm:counting}, with $\beta_d=\gamma_d+\frac d2\log\frac d{d-1}$ and the Gaussian integral from~\eqref{eq:rho}:
      \[
            \int_E\psi_E(\vu)^{\frac{d-1}{d}}\dd\vu
            =(2\pi)^{\frac{d-1}{2d}}
            \left(\frac d{d-1}\right)^{\frac{d-1}{2}}.
      \]
\end{remark}
\begin{remark}
      The following heuristic explains the exponent $p=\frac{d-1}{d}$ in~\eqref{eq:weighted}. We tile the floor $E$ with equal-sized tiles. For $r$ points in a single $(d-1)$-dimensional floor tile, the typical separation is of order $r^{-1/(d-1)}$. The antichain constraint restricts neighboring height differences to this scale, suggesting a volume factor $r^{-r/(d-1)}$. Label allocations contribute $1/r!$ per tile, of order $r^{-r}$ up to exponential factors. Writing $w_j$ for the approximately constant floor weight $\psi_E e^\sigma$ on tile $j$, the tile contributes roughly $w_j^r/r^{r/p}$. Maximizing the product over occupancies $r_j$ with $\sum_jr_j=n$ suggests $r_j\propto w_j^p$. Thus, asymptotically, conditioning on the antichain event replaces the original floor weights by their normalized $p$th powers. Sections~\ref{sec:counting_proof}--\ref{sec:lower-bound} make this heuristic precise.
\end{remark}

Section~\ref{sec:shortcuts} derives empirical concentration from Proposition~\ref{prop:volume}: varying $c$ controls heights, while varying $\sigma$ identifies the floor distribution. Together with the antichain condition, this concentration yields Theorem~\ref{thm:shape}.

Sections~\ref{sec:counting_proof}--\ref{sec:lower-bound} prove Proposition~\ref{prop:volume}. We tile $E$ by small reflected copies of a simplex $\Delta$. The one-tile volume estimates define $\gamma_d$ in Section~\ref{sec:counting_proof}. For the upper bound in Section~\ref{sec:upper-bound}, we discard comparisons between different tiles. For the lower bound in Section~\ref{sec:lower-bound}, we retain configurations away from tile boundaries that admit height extensions with common boundary values. Reflections match these values across adjacent tiles, allowing the independently chosen interiors to form one antichain.

\section{From weighted expectations to the limit shape}\label{sec:shortcuts}

We now identify the distribution of points along the limit shape and then deduce the Hausdorff statement. We assume Proposition~\ref{prop:volume}, whose proof occupies Sections~\ref{sec:counting_proof}--\ref{sec:lower-bound}.

\subsection{The limiting distribution}

Let $\vY\in\R^d$ have independent standard normal coordinates, and put
\begin{equation}\label{eq:Fstar_definition}
      \vF^\ast=\sqrt{\frac d{d-1}}\Proj_E\vY.
\end{equation}
Writing $\bar Y=\frac1d\sum_iY_i$ gives $F_i^\ast=\sqrt{\frac{d}{d-1}}(Y_i-\bar Y)$. Hence by direct calculation, every coordinate of $\vF^\ast$ has variance $1$, distinct coordinates have covariance $-\frac1{d-1}$, and the sum of the coordinates is zero. The density $\rho^\ast$ of $\vF^\ast$ on $E$ is
\begin{equation}\label{eq:rho-star}
      \rho^\ast(\vu)=\left(\frac{d-1}{2\pi d}\right)^{\frac{d-1}{2}}
      e^{-\frac{d-1}{2d}\|\vu\|^2},\qquad \vu\in E.
\end{equation}
Comparing with~\eqref{eq:rho}, we see that $\rho^\ast$ is the normalized power $\psi_E^{(d-1)/d}$.

The empirical measure of the sample and the limiting measure are
\begin{equation}\label{eq:mu_n_hat_definition}
      \widehat\mu_n=\frac1n\sum_{k=1}^n\delta_{\vX^{(k)}}, \qquad \mu^\ast=\law\bigl(\Phi(\vF^\ast)\bigr).
\end{equation}
Here $\Phi$ acts coordinatewise, and $\delta_{\vx}$ denotes unit mass at $\vx$, so $\widehat\mu_n(B)$ is the fraction of sampled points in $B$.
Each coordinate marginal of $\mu^\ast$ is uniform, because $\P(\Phi(Z)\le x)=x$ for a standard normal $Z$ and $0<x<1$. Recall that under $\P_n$, the sample is conditioned to form an antichain.

Let $\mathcal P_d$ denote the set of probability measures on $[0,1]^d$, equipped with the weak topology generated by the maps $\nu\mapsto\int g\dd\nu$ for continuous $g$; see~\cite[Chapter~1]{Billingsley}.

\begin{theorem}
      \label{prop:empirical-limit}\label{thm:empirical-limit}\label{thm:distribution}
      For every open neighborhood $\cU$ of $\mu^\ast$ in $\mathcal P_d$, there is $\alpha=\alpha(d,\cU)>0$ such that
      \begin{equation}\label{eq:empirical}
            \P_n(\widehat\mu_n\notin\cU)\le e^{-\alpha n}
      \end{equation}
      for all sufficiently large $n$. In particular, for every Borel set $B$ with $\mu^\ast(\partial B)=0$,
      \[
            \frac1n\#\{k:\vX^{(k)}\in B\}\longrightarrow\mu^\ast(B)
            \quad\text{in }\P_n\text{-probability}.
      \]
\end{theorem}
\begin{proof}
      We start by extracting conditioned expectations from Proposition~\ref{prop:volume}.
      For every bounded continuous $\sigma:E\to\R$ and $c<1/2$, applying~\eqref{eq:weighted} at $(\sigma,c)$ and $(0,0)$ and subtracting cancels the terms $\gamma_d-\frac12\log(2\pi)$. Since $\rho^\ast$ is the normalized power $\psi_E^{(d-1)/d}$, this gives
      \[
            \lim_{n\to\infty}\frac1n \log \frac{\cZ_n(\sigma,c)}{\cZ_n(0,0)} = \frac{d}{d-1}\log\frac{\int_E(\psi_E(\vu)e^{\sigma(\vu)})^{\frac{d-1}d}\dd\vu}{\int_E \psi_E(\vu)^{\frac{d-1}d}\dd\vu}=\frac{d}{d-1}\log\E e^{\frac{d-1}{d}\sigma(\vF^\ast)}.
      \]
      Write $\E_n=\E(\,\cdot\mid\cA_n)$ and $T_n=\sum_{k=1}^n(H^{(k)})^2$. Since $\cZ_n(0,0)=a_n$, definition~\eqref{eq:def-zn} identifies the ratio above with $\E_n\exp\{\sum_{k=1}^n\sigma(\vF^{(k)})+cT_n\}$. Therefore,
      \begin{equation}\label{eq:laplace}
            \lim_{n\to\infty}\frac1n\log\E_n \exp\biggl\{\sum_{k=1}^n\sigma(\vF^{(k)})+cT_n\biggr\}
            =\frac d{d-1}\log\E e^{\frac{d-1}{d}\sigma(\vF^\ast)}.
      \end{equation}
      This is the main tool that we will use to prove the theorem. We use~\eqref{eq:laplace} with $\sigma=0$ to control heights and with $c=0$ to control floor positions.

      Returning to the statement of the theorem, every open weak neighborhood $\cU$ of $\mu^\ast$ contains a basic open neighborhood $\cB$ of $\mu^\ast$ specified by the averages of finitely many continuous functions $g_1,\ldots,g_m$ on the cube:
      \[\cB=\left\{\mu\in\mathcal{P}_{d}:\int_{[0,1]^{d}}g_{j}\dd\mu-\int_{[0,1]^{d}}g_{j}\dd\mu^\ast<\eps_{j}\quad\text{for $j=1,\ldots,m$}\right\},\]%
      where $\eps_{1},\ldots,\eps_{m}$ are positive reals. By~\eqref{eq:mu_n_hat_definition}, we have
      \[\int_{[0,1]^{d}}g_{j}\dd\mu^\ast=\E g_{j}(\Phi(\vF^\ast))\qquad \text{and}\qquad \int_{[0,1]^{d}}g_{j}\dd\widehat\mu_n=\frac{1}{n}\sum_{k} g_{j}(\vX^{(k)}).\]
      Hence, we can write
      \[
            \P_n(\widehat\mu_n\notin\cU) \le \P_n(\widehat\mu_n\notin\cB) \le \sum_{j=1}^m \P_n\left(\frac{1}{n}\sum_{k} g_{j}(\vX^{(k)})-\E g_{j}(\Phi(\vF^\ast))\geq\eps_{j}\right).
      \]
      It now suffices to show that the probability on the right decreases exponentially for each continuous function $g_{j}$ and each real $\eps_{j}>0$. By the Stone--Weierstrass Theorem, the set of Lipschitz functions on $[0,1]^{d}$ is dense in the set of continuous functions. Hence, there is a Lipschitz function $\widehat{g}_{j}$ on $[0,1]^{d}$ such that $|g_{j}(\vx)-\widehat{g}_{j}(\vx)|<\frac{\eps_{j}}{4}$. Each of the empirical and limiting averages then changes by less than $\eps_j/4$, so
      \[\P_n\left(\frac{1}{n}\sum_{k} g_{j}(\vX^{(k)})-\E g_{j}(\Phi(\vF^\ast))\geq\eps_{j}\right)\leq \P_n\left(\frac{1}{n}\sum_{k} \widehat{g}_{j}(\vX^{(k)})-\E \widehat{g}_{j}(\Phi(\vF^\ast))\geq\frac{\eps_{j}}{2}\right).\]
      Write $\sigma_{j}=\widehat{g}_{j}\circ\Phi$ on $\R^d$.
      Recall that $\vX^{(k)}=\Phi(\vG^{(k)})=\Phi(\vF^{(k)}+H^{(k)}\ve).$ Then it suffices to prove that for each $j$, there is some $\alpha_{j}>0$ such that, for sufficiently large $n$:
      \[\P_n\left(\frac{1}{n}\sum_{k} \sigma_{j}(\vF^{(k)}+H^{(k)}\ve)-\E \sigma_{j}(\vF^\ast)\geq\frac{\eps_{j}}{2}\right)<e^{-\alpha_{j}n}.\]
      The function $\sigma_j$ is bounded and Lipschitz: $\widehat g_j$ has these properties, and $0<\Phi'\le(2\pi)^{-1/2}$ makes the coordinatewise map $\Phi$ Lipschitz as well.

      We now bound the mean squared height $T_n/n$ to control the error from replacing each point by its floor position. We fix $0<c<1/2$. Substituting $\sigma=0$ into~\eqref{eq:laplace} yields
      \[
            \log\E_n e^{cT_n}=o(n).
      \]
      In the event that $T_n\ge\delta n$, we have $e^{cT_n}\ge e^{c\delta n}$. Markov's inequality gives, for every fixed $\delta>0$,
      \begin{equation}\label{eq:height-tail}
            \P_n(T_n\ge\delta n)\le e^{-c\delta n}\E_n(e^{c T_n})= e^{-c\delta n+o(n)}.
      \end{equation}
      Now, recall that $\sigma_{j}$ is $L_{j}$-Lipschitz for some $L_{j}>0$.
      The Lipschitz bound gives
      \[\frac{1}{n}\sum_{k}\sigma_{j}(\vF^{(k)}+H^{(k)}\ve)- \frac{1}{n}\sum_{k}\sigma_{j}(\vF^{(k)})\leq L_{j}\frac{1}{n}\sum_{k}|H^{(k)}|\leq L_{j}\sqrt{\frac{T_n}{n}},\]
      where the second inequality is by Cauchy--Schwarz, $(\sum_k|H^{(k)}|)^2\le n\sum_k(H^{(k)})^2$. Consequently,
      \begin{equation}\label{eq:heights_neligible}\P_n\left(\frac{1}{n}\sum_{k}\sigma_{j}(\vF^{(k)}+H^{(k)}\ve)- \frac{1}{n}\sum_{k}\sigma_{j}(\vF^{(k)})\geq\frac{\eps_{j}}{4}\right)\leq \P_n\left(\frac{T_n}{n}\geq\left(\frac{\eps_{j}}{4 L_{j}}\right)^{2}\right)\le e^{-\beta_{j}n+o(n)},\end{equation}
      by~\eqref{eq:height-tail}, where $\beta_{j}=c\left(\frac{\eps_{j}}{4 L_{j}}\right)^{2}$. We next obtain a similar bound for the probability
      \[\P_n\left(\frac{1}{n}\sum_{k} \sigma_{j}(\vF^{(k)})-\E \sigma_{j}(\vF^\ast)\geq\frac{\eps_{j}}{4}\right),\]
      by analyzing the distribution of the floor values $\vF^{(k)}$. Set $\delta=\frac{\eps_{j}}{4}$ and define $\sigma:E\to \R$ by $\sigma(\vx)=\theta(\sigma_{j}(\vx)-\E \sigma_{j}(\vF^\ast))$ for some $\theta>0$, so $\E \sigma(\vF^\ast)=0$ and $\sigma$ is bounded and continuous.
      Markov's inequality gives, for every $\theta>0$,
      \[
            \P_n \left(\frac1n\sum_k  \sigma(\vF^{(k)})\geq\theta\delta\right)\le e^{-\theta\delta n} \E_n \left(\exp\left\{\sum_k \sigma(\vF^{(k)})\right\}\right).
      \]
      Combining this with~\eqref{eq:laplace} at $c=0$ then yields
      \begin{equation}\label{eq:deviations_sigmaF}
            \limsup_{n\to\infty}\frac{1}{n}\log\P_n \left(\frac1n\sum_k \sigma(\vF^{(k)})\geq\theta\delta\right)\le -\theta\delta+\frac d{d-1}\log\E e^{\frac{d-1}{d}\sigma(\vF^\ast)}.
      \end{equation}
      We now show that this bound is negative for a suitable $\theta>0$. Put $S_j=\sigma_j(\vF^\ast)-\E\sigma_j(\vF^\ast)$. This variable is bounded and has mean zero. Taylor expansion, uniform over its bounded range, gives $\E e^{tS_j}=1+t\E S_j+O(t^2)=1+O(t^2)$ and hence $\log\E e^{tS_j}=O(t^2)$ as $t\to0$. Since $\sigma(\vF^\ast)=\theta S_j$, the right-hand side of~\eqref{eq:deviations_sigmaF} is $-\theta\delta+O(\theta^2)$, which is negative for sufficiently small $\theta>0$, independently of $n$. Thus, for some $\gamma_j>0$,
      \[
            \P_n\left(\frac1n\sum_k\sigma_j(\vF^{(k)})-\E\sigma_j(\vF^\ast)\ge\frac{\eps_j}{4}\right)\le e^{-\gamma_jn+o(n)}.
      \]
      Combining this with the height bound in~\eqref{eq:heights_neligible} gives
      \[\P_n\left(\frac{1}{n}\sum_{k} \sigma_{j}(\vF^{(k)}+H^{(k)}\ve)-\E \sigma_{j}(\vF^\ast)\geq\frac{\eps_{j}}{2}\right)\leq e^{-\beta_{j}n+o(n)}+e^{-\gamma_{j}n+o(n)},\]
      which decreases exponentially in $n$. Summing over the finitely many tests gives
      \[
            \P_n(\widehat\mu_n\notin\cU) \le\sum_{j=1}^m \bigl(e^{-\beta_{j}n+o(n)}+e^{-\gamma_{j}n+o(n)}\bigr)\le e^{-\alpha n}
      \]
      for some $\alpha>0$ and all sufficiently large $n$. Finally, the second assertion follows from~\cite[Theorem~2.1]{Billingsley}, which gives continuity of $\nu\mapsto\nu(B)$ at $\mu^\ast$ whenever $\mu^\ast(\partial B)=0$.
\end{proof}

\subsection{From empirical concentration to the shape}

In this section, we prove Theorem~\ref{thm:shape} using Theorem~\ref{thm:distribution}. To exclude isolated outliers, we use the following comparison property of the support.
\begin{lemma}\label{lem:support}
      The support of $\mu^\ast$ is
      \begin{equation}\label{eq:closure}
            \overline{\cS_d} =\cS_d\cup\{\vx\in[0,1]^d:\min_i\vx_i=0,\ \max_i\vx_i=1\}.
      \end{equation}
      For every $\vx\in[0,1]^d\setminus\overline{\cS_d}$, there is $\vy\in\cS_d$ such that $\vx\prec\vy$ or $\vy\prec\vx$.
\end{lemma}
\begin{proof}
      Since $\rho^\ast$ is positive on $E$, the support is $\overline{\Phi(E)}=\overline{\cS_d}$. Continuity in the open cube and the endpoint limits of $\Phi^{-1}$ give~\eqref{eq:closure}: a boundary limit must have both a zero and a one, since otherwise the sum of the inverse-normal coordinates diverges with one sign. Conversely, let $I_b=\{i:\vx_i=b\}$ for $b\in\{0,1\}$ be nonempty, and put $C=\sum_{i\notin I_0\cup I_1}\Phi^{-1}(\vx_i)$. Assign inverse-normal coordinates $-t/|I_0|$ on $I_0$ and $(t-C)/|I_1|$ on $I_1$, leaving the others fixed. Their sum is zero, and their $\Phi$-images approach the required boundary point as $t\to\infty$.

      For $\vx\notin\overline{\cS_d}$, we have $\max_i\vx_i-\min_i\vx_i<1$. The strictly increasing continuous function
      \[
            t\longmapsto\sum_{i=1}^d\Phi^{-1}(\vx_i+t), \qquad -\min_i\vx_i<t<1-\max_i\vx_i,
      \]
      tends to $-\infty$ and $+\infty$ at the two endpoints. Therefore, it vanishes at some value $t_0$ which is nonzero as otherwise $\vx\in\cS_d$. The point $\vy=\vx+t_0(1,\ldots,1)$ lies in $\cS_d$ and is strictly comparable with $\vx$.
\end{proof}

\begin{proof}[Proof of Theorem~\ref{thm:shape}]
      Theorem~\ref{thm:empirical-limit} gives empirical concentration, so it remains to obtain the corresponding Hausdorff bound. Fix $\delta>0$. We construct finitely many continuous test functions $g_{1},\ldots,g_{m}$ on $[0,1]^{d}$ that detect any hole or outlier at distance greater than $\delta$. More precisely, these test functions all have positive $\mu^\ast$-integrals but, if there is a hole or outlier at distance greater than $\delta$, then at least one of these test functions will have empirical integral $0$.

      For holes, cover $\overline{\cS_d}$ by finitely many balls of radius $\frac\delta2$ centered at $\vy^{(j)}\in\overline{\cS_d}$. The functions
      \[
            g_j(\vx)=\left(\frac\delta2-\|\vx-\vy^{(j)}\|\right) \vee 0
      \]
      have positive $\mu^\ast$-integrals by Lemma~\ref{lem:support}. If some $\vy\in\overline{\cS_d}$ has distance at least $\delta$ from the sample, a center within $\frac\delta2$ of $\vy$ is at distance at least $\frac\delta2$ from every sample point. The corresponding $g_j$ then has empirical integral zero.

      For outliers, the set
      \[
            K_\delta=\{\vx\in[0,1]^d:\dist(\vx,\overline{\cS_d})\ge\delta\}
      \]
      is compact. For each $\vx\in K_\delta$, Lemma~\ref{lem:support} supplies a strictly comparable point $\vy_{\vx}\in\cS_d$. Strict inequalities persist in neighborhoods, so there are neighborhoods $U_{\vx}$ and $V_{\vx}$ of these points such that every point in $U_{\vx}$ is strictly comparable with every point in $V_{\vx}$. Since $K_\delta$ is compact, we may choose a finite subcover $(U_j)$ with corresponding neighborhoods $(V_j)$. As above, choose a nonnegative continuous test function supported in each $V_j$, positive at the chosen point of $\cS_d$. Its $\mu^\ast$-integral is positive, whereas its empirical integral is zero for every antichain meeting $U_j$.

      Relabel the previous test functions as $g_1,\ldots,g_m$ and let
      \[
            \cU_\delta=\left\{\nu: \int g_j\dd\nu>\frac12\int g_j\dd\mu^\ast \text{ for }1\le j\le m\right\}.
      \]
      This is an open weak neighborhood of $\mu^\ast$, fixed independently of $n$. A hole or outlier makes one of these averages zero, forcing $\widehat\mu_n\notin\cU_\delta$. Theorem~\ref{thm:empirical-limit} therefore gives
      \[
            \P_n\bigl(d_{\rm H}(A_n,\overline{\cS_d})>\delta\bigr) \le\P_n(\widehat\mu_n\notin\cU_\delta) \le e^{-\alpha n}
      \]
      for some $\alpha=\alpha(d,\delta)>0$ and all sufficiently large $n$.
\end{proof}

\section{Asymptotic results on tiles}\label{sec:counting_proof}

We develop the local geometry needed to prove the weighted estimate in Proposition~\ref{prop:volume}. First we describe the constraints on the heights $H^{(k)}$ in an antichain in terms of the floor positions $\vF^{(k)}$. We then describe a tiling of the floor hyperplane $E$ with simplices. Finally we determine (unweighted) volume estimates within one tile after fixing one height. These estimates define the finite constant $\gamma_d$ in Proposition~\ref{prop:volume}. The weights $\sigma$ and $c$ enter in Section~\ref{sec:upper-bound}.
\subsection{Reinterpretation of antichains}
Recall that $\ve=\frac{1}{\sqrt d}(1,\ldots,1)$ and $E=\ve^\perp$.
Before describing the tiling, we describe antichains in terms of the coordinates $\vf^{(k)}$ and $h^{(k)}$ denoting floor positions and heights respectively. To do this we define the function $\mathcal L:E\to[0,\infty)$ as follows:
\[
      \Ldir{\vu}=-\sqrt d\min_{1\le i\le d}\vu_i, \qquad  \vu\in E.
\]
The function $\mathcal L$ gives an upper bound on the allowed height increment for a given floor displacement.

\begin{proposition} \label{prop:L_definition}
      For $\ff=(\vf^{(1)},\ldots,\vf^{(r)})\in E^r$ and $\fh=(h^{(1)},\ldots,h^{(r)})\in\R^r$, the points $\vf^{(k)}+h^{(k)}\ve$, $1\le k\le r$, form an antichain if and only if
      \begin{equation}\label{eq:constraints}
            h^{(k)}-h^{(\ell)}\le\Ldir{\vf^{(k)}-\vf^{(\ell)}} \qquad\text{for every ordered pair }(k,\ell).
      \end{equation}
\end{proposition}
\begin{proof}
      For an ordered pair $(k,\ell)$, we have $\vf^{(k)}+h^{(k)}\ve \succ \vf^{(\ell)}+h^{(\ell)}\ve$ exactly when
      \[
            h^{(k)}-h^{(\ell)}>
            -\sqrt d\min_i(\vf^{(k)}-\vf^{(\ell)})_i
            =\Ldir{\vf^{(k)}-\vf^{(\ell)}}.
      \]
      Excluding this comparison for every ordered pair gives~\eqref{eq:constraints}.
\end{proof}

We now prove some basic properties of $\mathcal L$.
\begin{lemma}
      For $\vu,\vv\in E$ and $s\ge0$,
      \begin{equation}\label{eq:height-bound}
            \Ldir{\vu+\vv}\le\Ldir{\vu}+\Ldir{\vv},\qquad \Ldir{s\vu}=s\Ldir{\vu},\qquad \frac{\|\vu\|}{\sqrt{d-1}}\le\Ldir{\vu} \le\sqrt{d-1}\|\vu\|.
      \end{equation}
      The function $\mathcal L$ is invariant under coordinate permutations.
\end{lemma}
\begin{proof}
      Subadditivity, positive homogeneity, and permutation invariance follow from the definition. The coordinates $\vu_i$ sum to $0$ and each term $\vu_i+\frac{\Ldir{\vu}}{\sqrt{d}}$ is nonnegative, so
      \[
            \|\vu\|^2 =\sum_{i=1}^d\Bigl(\vu_i+\frac{\Ldir{\vu}}{\sqrt{d}}\Bigr)^2-\Ldir{\vu}^2 \le \biggl(\sum_{i=1}^d\Bigl(\vu_i+\frac{\Ldir{\vu}}{\sqrt{d}}\Bigr)\biggr)^2-\Ldir{\vu}^2 = (d-1)\Ldir{\vu}^2.
      \]
      This proves the lower bound. For the upper bound, the Cauchy--Schwarz inequality implies
      \[
            \vu_i^2=\Bigl(\sum_{j\ne i}\vu_j\Bigr)^2 \le(d-1)\sum_{j\ne i}\vu_j^2,
      \]
      and thus $\Ldir{\vu}\le\sqrt{d}\max_{1\le i\le d}|\vu_i|\le\sqrt{d-1}\|\vu\|$.
\end{proof}

Applying~\eqref{eq:constraints} in both orders and using~\eqref{eq:height-bound} gives the following Lipschitz property.

\begin{lemma}\label{lem:height}\label{cor:height}
      If the points $\vf^{(k)}+h^{(k)}\ve$, $1\le k\le r$, form an antichain then they satisfy the Lipschitz condition $|h^{(k)}-h^{(\ell)}|\le\sqrt{d-1}\|\vf^{(k)}-\vf^{(\ell)}\|$.
\end{lemma}

\subsection{Tiling of \texorpdfstring{$E$}{E} with simplices}\label{sec:tiling}
\begin{figure}[!htbp]
      \centering
      \resizebox{0.60\linewidth}{!}{%
            \begin{tikzpicture}[
                        x=3.55cm, y=3.55cm,
                        every node/.style={font=\small,inner sep=1pt},
                        tiling wall/.style={dash pattern=on 3.2pt off 3.2pt,line width=0.45pt},
                        side label/.style={fill=white,inner sep=1.2pt,font=\footnotesize}
                  ]
                  \pgfmathsetmacro{\tileH}{sqrt(3)/2}
                  \path[use as bounding box] (-1.08,{-0.43*\tileH}) rectangle (2.58,{2.17*\tileH});
                  \begin{scope}
                        \clip (-1.08,{-0.43*\tileH}) rectangle (2.58,{2.17*\tileH});
                        \fill[yellow!12] (0,0)--(1,0)--(0.5,\tileH)--cycle;
                        \fill[black!18] (1,{2*\tileH})--(2,{2*\tileH})--(1.5,\tileH)--cycle;
                        \foreach \j in {-4,...,5} {
                                    \draw[tiling wall] (-4,{\j*\tileH})--(5,{\j*\tileH});
                              }
                        \foreach \j in {-5,...,5} {
                                    \draw[tiling wall] ({\j-4},{-4*sqrt(3)})--({\j+4},{4*sqrt(3)});
                                    \draw[tiling wall] ({\j-4},{4*sqrt(3)})--({\j+4},{-4*sqrt(3)});
                              }
                        \draw[line width=0.8pt] (0,0)--(1,0)--(0.5,\tileH)--cycle;
                        \draw[line width=0.8pt] (1,{2*\tileH})--(2,{2*\tileH})--(1.5,\tileH)--cycle;

                        \node[font=\large] at (0.5,{\tileH/3}) {$\Delta$};
                        \node at (0,{2*\tileH/3}) {$r_1\Delta$};
                        \node at (1,{2*\tileH/3}) {$r_0\Delta$};
                        \node at (-0.5,{\tileH/3}) {$r_1r_2\Delta$};
                        \node at (1.5,{\tileH/3}) {$r_0r_2\Delta$};
                        \node at (2,{2*\tileH/3}) {$r_0r_2r_1\Delta$};
                        \node at (0,{4*\tileH/3}) {$r_1r_0\Delta$};
                        \node at (1,{4*\tileH/3}) {$r_0r_1\Delta$};
                        \node at (2,{4*\tileH/3}) {$r_0r_1r_2r_1\Delta$};
                        \node at (-0.5,{5*\tileH/3}) {$r_1r_0r_2\Delta$};
                        \node at (0.5,{5*\tileH/3}) {$r_0r_1r_0\Delta$};
                        \node[align=center] at (1.5,{5*\tileH/3}) {\large $\Delta'$\\[-1pt]$=r_0r_1r_2\Delta$};

                        \node at (-0.5,{-0.25*\tileH}) {$r_1r_2r_1\Delta$};
                        \node at (0.5,{-0.25*\tileH}) {$r_2\Delta$};
                        \node at (1.5,{-0.25*\tileH}) {$r_0r_2r_0\Delta$};
                        \node[side label] at (0.25,{\tileH/2}) {$r_1$};
                        \node[side label] at (0.75,{\tileH/2}) {$r_0$};
                        \node[side label] at (0.5,0) {$r_2$};
                  \end{scope}
            \end{tikzpicture}}
      \caption{A view of the tiling of the floor $E=\{(x,y,z):x+y+z=0\}$ when $d=3$. The yellow triangle is the reference tile $\Delta$, and the gray triangle is $\Delta'=r_0r_1r_2\Delta$. The generators $r_i$ are the fixed reflections in the three sides of $\Delta$ and products act from right to left.}
      \label{fig:triangles}
\end{figure}
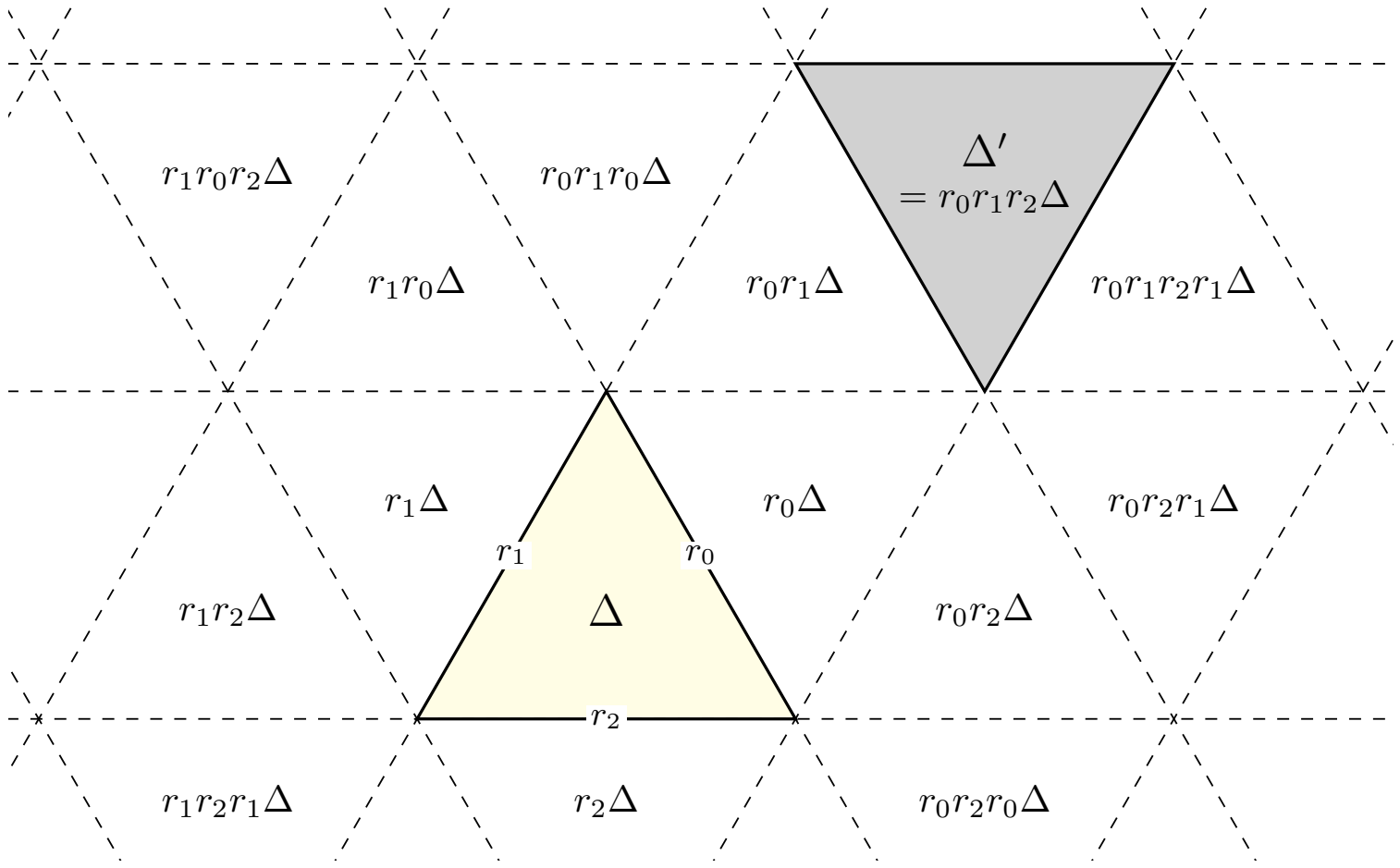

We will use a tiling of the space $E$ by simplices. This tiling is the Coxeter complex for the affine Coxeter group $\widetilde{A}_{d-1}$. Consider the reference simplex
\[
      \Delta=\{\vf\in E:\vf_1\ge\cdots\ge\vf_d,\quad \vf_1-\vf_d\le1\}.
\]
Its vertices are $\mathbf0$ and, for $1\le j<d$, the vectors with first $j$ coordinates $\frac{d-j}{d}$ and remaining coordinates $-\frac jd$. For $d=3$ it is an equilateral triangle, as depicted in Figure~\ref{fig:triangles}.
Consider the family $\cT$ of images of $\Delta$ under compositions of its side reflections. Its tiles $\Delta'$ have disjoint interiors and are translates of coordinate permutations of $\Delta$; see Proposition~\ref{prop:tiling}.
For $1\le i<d$, let $r_i$ be reflection in $\vf_i=\vf_{i+1}$, which exchanges those coordinates. The remaining side reflection of $\Delta$ is
\[
      r_0(\vf)=(\vf_d+1,\vf_2,\ldots,\vf_{d-1},\vf_1-1).
\]
Importantly, these reflections preserve $\mathcal L$, in the sense that for $\vf^{(k)},\vf^{(\ell)}\in E$, we have
\[\Ldir{r_{i}(\vf^{(k)})-r_{i}(\vf^{(\ell)})}=\Ldir{\vf^{(k)}-\vf^{(\ell)}},\]
for $0\leq i\leq d-1$. The translation part cancels in the difference, leaving a coordinate permutation. Figure~\ref{fig:triangles} shows part of the tiling in dimension three, with a nonadjacent tile $\Delta'$ highlighted.

For a measurable set $B\subset E$, write $|B|$ for its $(d-1)$-dimensional Euclidean volume. Denote by $\Omega=\langle r_0,\ldots,r_{d-1}\rangle$ the group of affine isometries of $E$ generated by the side reflections above.
\begin{proposition}\label{prop:tiling}
      The reflected copies of $\Delta$ tile $E$ with disjoint interiors, and every compact set meets only finitely many tiles. For each tile $\Delta'$, there is a unique map $\omega\in \Omega$ generated by the side reflections such that $\Delta'=\omega\Delta$. For each $x\in E$ there is a unique point $y\in\Delta$ such that $\omega y=x$ for some $\omega\in \Omega$. For any fixed cube $\square\subset E$ in orthonormal coordinates, the number of scale-$\eps$ tiles $\eps\Delta'$ completely contained in $\square$ satisfies
      \begin{equation}\label{eq:tile-number}
            q(\eps,\square)=\frac{|\square|}{\eps^{d-1}|\Delta|}+O(\eps^{2-d})\qquad(\eps\downarrow0).
      \end{equation}
\end{proposition}
\begin{proof} In our coordinates, the affine reflection hyperplanes of type $\widetilde{A}_{d-1}$ are
      \[
            H_{ij,k}=\{\vf\in E:\vf_i-\vf_j=k\},\qquad 1\le i<j\le d,\quad k\in\mathbb Z.
      \]
      These are the hyperplanes $H_{\alpha,k}$ in~\cite[Section~4.1]{Humphreys} for the roots $\alpha=\mathbf e_i-\mathbf e_j$, where $\mathbf e_1,\ldots,\mathbf e_d$ are the standard basis vectors of $\R^d$; see~\cite[Section~2.10]{Humphreys}. The connected components of their complement are called alcoves. The interior of $\Delta$ is precisely the region where $0<\vf_i-\vf_j<1$ for every $i<j$, which is the fundamental alcove described in~\cite[Section~4.3]{Humphreys}. Its wall reflections are $r_0,\ldots,r_{d-1}$. By~\cite[Proposition~4.3]{Humphreys}, the affine Weyl group is generated by these reflections and hence equals $\Omega$. Thus the tiling and uniqueness of the maps $\omega$ follow from~\cite[Proposition~4.3 and Theorem~4.5(c)]{Humphreys}. Also, it follows from~\cite[Theorem~4.8]{Humphreys} that the closed alcove $\Delta$ meets each group orbit in exactly one point. Finally, by~\cite[Proposition~4.2]{Humphreys}, each tile is a translate of a coordinate permutation of $\Delta$ by a vector in $E\cap\mathbb Z^d$. Since $E\cap\mathbb Z^d$ is discrete, $\Delta$ is bounded, and there are only finitely many coordinate permutations, only finitely many such translates can meet a fixed compact set.

      The part of $\square$ not covered by complete tiles lies within distance $\eps\diam(\Delta)$ of its boundary, so its volume is $O(\eps)$. Thus for some $C_\square>0$,
      \[
            0\le |\square|-q(\eps,\square)\eps^{d-1}|\Delta| \le C_{\square}\eps.
      \]
      Dividing by $\eps^{d-1}|\Delta|$ proves~\eqref{eq:tile-number}.
\end{proof}

\subsection{One-tile estimates and exponential rates}
Let $B\subset E$ be measurable with $|B|>0$. We will estimate volumes of sets of antichains with floor values in $B$, then specialize to $B=\eps\Delta'$ for a tile $\Delta'\in\cT$.
For $r$ labeled points, write $\ff=(\vf^{(1)},\ldots,\vf^{(r)})\in E^r$ and $\fh=(h^{(1)},\ldots,h^{(r)})\in\R^r$. The configuration space just described is
\begin{equation}\label{eq:A_n_of_B_definition}
      \cA_r(B)=\{(\ff,\fh)\in B^r\times\R^r:
      \vf^{(k)}+h^{(k)}\ve\not\prec\vf^{(\ell)}+h^{(\ell)}\ve\text{ for all }k\ne\ell\}.
\end{equation}
We use product Lebesgue measure $\dd\ff\dd\fh=\prod_k\dd\vf^{(k)}\dd h^{(k)}$, with the induced Euclidean measure on $E$, and write $\Vol_q$ for $q$-dimensional volume.

The volume of $\cA_r(B)$ is infinite because translating all heights preserves admissibility. We remove this common translation by fixing the first height.
\begin{equation}\label{eq:pinned-set}
      \cA_r^0(B):=\{(\ff,\fh)\in\cA_r(B):h^{(1)}=0\},\qquad r\ge1.
\end{equation}
For fixed floor positions $\ff\in E^r$, we define the set $\cH_r^0(\ff)$ of heights forming antichains with these floor positions:
\[
      \cH_r^0(\ff):= \left\{\fh\in\R^r:(\ff,\fh)\in\cA_r(E),\ h^{(1)}=0\right\}.
\]
This set carries $(r-1)$-dimensional Lebesgue measure in the free coordinates $(h^{(2)},\ldots,h^{(r)})$. For fixed $\ff$, every admissible height list is uniquely obtained by the map
\[
      (s,\fh)\longmapsto (s,s+h^{(2)},\ldots,s+h^{(r)}), \qquad s\in\R,\quad\fh\in\cH_r^0(\ff).
\]
The Jacobian of this map in the free coordinates $(s,h^{(2)},\ldots,h^{(r)})$ is one. Define the pinned volume
\begin{equation}\label{eq:slice}
      V_r(B):=\Vol_{dr-1}\bigl(\cA_r^0(B)\bigr)
      =\int_{B^r} \Vol_{r-1}\bigl(\cH_r^0(\ff)\bigr)\dd\ff.
\end{equation}
In particular, $V_1(B)=|B|$. We apply this definition to $B=\eps\Delta'$ for a fixed tile $\Delta'\in \cT$ and $\eps>0$, whose volume is $|B|=\eps^{d-1}|\Delta|$.

\begin{lemma}
      For every integer $r\ge1$, $\eps>0$, and $\Delta'\in\cT$,
      \begin{equation}\label{eq:cut-scale}
            V_r(\eps\Delta') =\eps^{dr-1}V_r(\Delta') =\eps^{dr-1}V_r(\Delta).
      \end{equation}
\end{lemma}
\begin{proof}
      Translations and coordinate permutations preserve~\eqref{eq:constraints} and the Lebesgue measure. Scaling the floor and the $r-1$ heights $h^{(2)},\ldots,h^{(r)}$ by $\eps$ preserves the pin $h^{(1)}=0$ and has Jacobian $\eps^{(d-1)r}\eps^{r-1}=\eps^{dr-1}$, proving~\eqref{eq:cut-scale}.
\end{proof}

The next lemma gives $\log V_r(\Delta)=-\frac r{d-1}\log r+O_d(r)$ and defines $\gamma_d$ using a limit superior. We give a direct proof, while the analogous cube estimate appears in~\cite{GP}.
\begin{lemma}\label{lem:constant}
      For $r\in\N=\{1,2,\ldots\}$, set
      \begin{equation}\label{eq:gamma}
            \gamma_{d,r}=\frac1r\log V_r(\Delta)+\frac1{d-1}\log\frac r{|\Delta|^d}, \qquad \gamma_d=\limsup_{r\to\infty}\gamma_{d,r}.
      \end{equation}
      The sequence $(\gamma_{d,r})_{r\in\N}$ is uniformly bounded, so $\gamma_d$ is finite.
\end{lemma}
\begin{proof}
      For the upper bound, fix distinct floor positions in $\Delta$. The quantity $\|\vu\|_{\mathcal L}:=\Ldir{\vu}+\Ldir{-\vu}=\sqrt d(\max_i\vu_i-\min_i\vu_i)$ is a norm on $E$ and is equal to the width of the allowed height-difference interval for floor displacement $\vu$ (by~\eqref{eq:constraints}, applied in both orders). Write $B_{\mathcal L}(\vu,s)$ for the open ball in the ${\mathcal L}$-norm. Choose a spanning tree on the floor positions that minimizes the total edge length, with edge lengths $L_1\ge\cdots\ge L_{r-1}$ in the ${\mathcal L}$-norm.
      Root the tree at vertex $1$ and fix $h^{(1)}=0$. After discarding all non-tree constraints, each edge of length $L_j$ restricts the corresponding child's height to an interval of length $L_j$ for any fixed height of its parent. Integrating the heights from the leaves toward the root therefore gives the upper bound $\prod_jL_j$ for the height measure.
      We now bound each edge length $L_j$ in terms of $R$, where $R=\max_{\vu,\vv\in\Delta}\|\vu-\vv\|_{\mathcal L}$ is the diameter of $\Delta$.
      Remove the edges $L_1,\ldots,L_{j-1},L_{j}$, splitting the tree into $j+1$ components. Then any two vertices in different remaining components have distance at least $L_{j}$. Otherwise, adding a shorter cross-component edge to the original tree, then deleting an edge in $L_1,\ldots,L_{j-1},L_{j}$ from the resulting cycle would shorten the tree.

      Choose one vertex from each component. The balls of radius $L_{j}/2$ centered at these representatives are disjoint and lie in $B_{\mathcal L}(\mathbf0,R+L_{j}/2)$, since $\mathbf0\in\Delta$. Comparing ball volumes, which scale as the $(d-1)$st power of the radius, gives
      \[
            (j+1)\left(L_j/2\right)^{d-1} \le\left(R+L_j/2\right)^{d-1} \le\left(3R/2\right)^{d-1}.
      \]
      Consequently, for some constant $C=C(d)>0$,
      \[
            L_j\le Cj^{-\frac1{d-1}},\qquad \prod_{j=1}^{r-1}L_j\le \frac{C^{r-1}}{((r-1)!)^{\frac1{d-1}}}.
      \]
      The bound is uniform over distinct floor positions, and coincidences between floor positions form a null set, so~\eqref{eq:slice} yields
      \[
            V_r(\Delta) \le\frac{|\Delta|^rC^{r-1}}{((r-1)!)^{\frac1{d-1}}}.
      \]
      Substituting this bound into~\eqref{eq:gamma} yields
      \[
            \gamma_{d,r}\le \frac{(r-1)\log C}{r} +\frac{r\log r-\log((r-1)!)}{(d-1)r} -\frac{\log|\Delta|}{d-1},
      \]
      which is uniformly bounded by Stirling's formula $\log(r!)=r\log r-r+O(\log(r+1))$.

      For the lower bound, choose Euclidean balls $B_1,\ldots,B_r\subset\Delta$ such that
      \[
            |B_i|\ge\frac{C|\Delta|}{r},\qquad \dist(B_i,B_j)\ge C r^{-\frac1{d-1}} \quad\text{ for }i\ne j
      \]
      for some $C=C(d)\in(0,1]$. To construct them, choose a fixed cube $\square$ inside $\Delta$, and let $C_{1}>0$ be its sidelength. Divide $\square$ into $m^{d-1}$ subcubes of sidelength $\frac{C_{1}}{m}$, where $m=\lceil r^{1/(d-1)}\rceil\le2r^{1/(d-1)}$. In any $r$ of these subcubes, take the ball of radius $\frac{C_{1}}{4m}$ which is concentric with the subcube. This gives both bounds with a fixed $C$. Requiring one floor point in each ball, with arbitrary labels, has volume at least $|\Delta|^rr!(C/r)^r$. By~\eqref{eq:height-bound}, the free heights $h^{(2)},\ldots,h^{(r)}$ in a common interval centered at zero of length $\frac C{\sqrt{d-1}}r^{-\frac1{d-1}}$ satisfy every height inequality with $h^{(1)}=0$. Therefore
      \[
            V_r(\Delta) \ge |\Delta|^rr!\left(\frac Cr\right)^r \left(\frac C{\sqrt{d-1}}r^{-\frac1{d-1}}\right)^{r-1}.
      \]
      Taking logarithms and using Stirling's formula gives $\log V_r(\Delta)\ge-\frac r{d-1}\log r-O_d(r)$. Substituting into~\eqref{eq:gamma} shows that $\gamma_{d,r}$ is uniformly bounded below.
\end{proof}

\section{The upper bound for Proposition~\ref{prop:volume}}

\label{sec:upper-bound}
We start by translating Proposition~\ref{prop:volume} to a more convenient form. Throughout this section, we fix a bounded continuous function $\sigma:E\to\R$ and $c<1/2$. Define the floor weight
\begin{equation}\label{eq:floor_weight_def}
      w(\vu)=\psi_E(\vu)e^{\sigma(\vu)}=(2\pi)^{-\frac{d-1}{2}}e^{\sigma(\vu)-\frac{\|\vu\|^2}{2}},\qquad \vu\in E.
\end{equation}
As described after~\eqref{eq:split}, the floors $\vF^{(k)}$ and heights $H^{(k)}$ follow independent standard Gaussian distributions (before conditioning). Hence, the weighted volume $\cZ_n(\sigma,c)$ defined in~\eqref{eq:def-zn} can be written as
\begin{equation}\label{eq:Z-integral}
      \cZ_n(\sigma,c)=(2\pi)^{-n/2}\int_{\cA_n(E)}
      e^{-(\frac12-c)\|\fh\|^2}\prod_{k=1}^n w(\vf^{(k)})\dd\ff\dd\fh.
\end{equation}
With this notation, Proposition~\ref{prop:volume} becomes
\begin{equation}\label{eq:integral-target}
      \lim_{n\to\infty}\left(\frac1n\log\cZ_n(\sigma,c)+\frac{\log n}{d-1}\right)
      =\gamma_d-\frac12\log(2\pi)
      +\frac d{d-1}\log\int_Ew(\vu)^{\frac{d-1}{d}}\dd\vu.
\end{equation}

The purpose of this section is to prove the upper bound for~\eqref{eq:integral-target}, namely \eqref{eq:upper-bound}. We do this by considering sets of points for which the restriction to any single tile is an antichain. To determine the corresponding weighted volume, we start with estimates in one floor tile, then we extend this to $E$ using a Poisson point process whose restrictions to different floor tiles are independent. Section~\ref{sec:lower-bound} proves the matching lower bound.

\subsection{Weighted contribution of one tile}
For a tile $B=\eps\Delta'$, set $w_B=\sup_Bw$. We estimate the contribution in~\eqref{eq:Z-integral} from $r\ge1$ labeled points with floors in $B$. The pinned volume $V_r(B)$ already integrates over their floor positions and relative heights; the next lemma integrates the remaining common height translation.

\begin{lemma}\label{lem:tile-height}
      For every $r\in\N$, $\eps>0$, $c<1/2$, and $\Delta'\in\cT$, put $B=\eps\Delta'$. Then
      \begin{equation}\label{eq:tile-weight-bound}
            \frac1{(2\pi)^{r/2}}\int_{\cA_r(B)}e^{-(\frac12-c)\|\fh\|^2}
            \prod_{k=1}^rw(\vf^{(k)})\dd\ff\dd\fh
            \le\frac{w_B^r}{(2\pi)^{r/2}}
            \int_{\cA_r(B)}e^{-(\frac12-c)\|\fh\|^2}\dd\ff\dd\fh
            \le w_B^rv_r(\eps,c),
      \end{equation}
      where $v_0(\eps,c):=1$ and, for $r\ge1$,
      \begin{equation}\label{eq:coefficients}
            v_r(\eps,c):=\frac{(2\pi)^{-\frac{r-1}{2}}}{\sqrt{r(1-2c)}}\eps^{dr-1} V_r(\Delta).
      \end{equation}
\end{lemma}
\begin{proof}
      The first inequality follows from $0\leq w(\vf^{(k)})\leq w_B$ for $\vf^{(k)}\in B$. For the second, we will use the parameterization of admissible heights preceding~\eqref{eq:slice}.
      Let $g(u)=\phi(u)e^{cu^2}$, so that
      $(2\pi)^{-r/2}e^{-(\frac12-c)\|\fh\|^2}
            =\prod_{k=1}^r g(h^{(k)})$. Then we can write
      \[\frac{1}{(2\pi)^{r/2}}
            \int_{\cA_r(B)}e^{-(\frac12-c)\|\fh\|^2}\dd\ff\dd\fh=
            \int_{B^{r}}\int_{\cH_r^0(\ff)}\int_{\R}\prod_{k=1}^r g(s+h^{(k)})\dd s\dd h^{(2)}\cdots \dd h^{(r)}\dd\ff.\]
      H\"older's inequality for $r$ functions with all exponents equal to $r$ and translation invariance of the resulting integral yield
      \[
            \int_\R\prod_{k=1}^r g(s+h^{(k)})\dd s
            \le \left(\prod_{k=1}^r\int_\R g(s+h^{(k)})^r\dd s\right)^{\frac{1}{r}}
            = \int_\R g(s)^r\dd s
            =\frac{(2\pi)^{-\frac{r-1}{2}}}
            {\sqrt{r(1-2c)}}.
      \]
      Integrating over $\fh\in\cH_r^0(\ff)$ and $\ff\in B^{r}=(\eps\Delta')^r$ multiplies this bound by
      \[
            \int_{(\eps\Delta')^r} \Vol_{r-1}\bigl(\cH_r^0(\ff)\bigr)\dd\ff
            =V_r(\eps\Delta')=\eps^{dr-1}V_r(\Delta),
      \]
      by~\eqref{eq:slice} and~\eqref{eq:cut-scale}. This gives the claimed bound.
\end{proof}

In particular, the values $v_r(\eps,c)$ are independent of the tile. The following lemma records their logarithmic asymptotics.
\begin{lemma}
      For fixed $\eps>0$ and $c<1/2$,
      \begin{equation}\label{eq:coefficient-rate}
            \limsup_{r\to\infty}\left(\frac{\log v_r(\eps,c)}{r}+\frac{\log r}{d-1}\right) =b_\eps:=\gamma_d-\frac12\log(2\pi)+\frac d{d-1}\log|\eps\Delta|.
      \end{equation}
\end{lemma}
\begin{proof}
      Substitute~\eqref{eq:gamma} into~\eqref{eq:coefficients}. Terms involving $\log r/r$, $\log\eps/r$, and $\log(1-2c)/r$ vanish. The remaining floor-volume term is $\frac d{d-1}\log|\Delta|+d\log\eps=\frac d{d-1}\log|\eps\Delta|$, because $|\varepsilon\Delta|=\varepsilon^{d-1}|\Delta|$.
\end{proof}

\subsection{Poissonization and independent tiles}

Fix $0<\eps\le1$. We use the Poisson point process $\Pi_t$ on $E\times\R$ with intensity
\[
      t\,w(\vf)\frac{e^{-(\frac12-c)h^2}}{\sqrt{2\pi}}\dd\vf\dd h,
      \qquad t\ge0.
\]
This means that counts in disjoint regions are independent Poisson random variables with means given by the intensity integral. Its total mean is $t\Lambda$, where
\begin{equation}\label{eq:poisson-mass}
      \Lambda=\frac1{\sqrt{1-2c}}\int_Ew(\vu)\dd\vu<\infty.
\end{equation}
Thus $\Pi_t$ contains finitely many points almost surely. We call it an antichain when the corresponding points $\vf+h\ve$ form an antichain; the empty configuration is allowed. Define
\[
      \mathcal Z(t)=e^{t\Lambda}\P(\Pi_t\text{ is an antichain}),
      \qquad \mathcal V(t)=\sum_{r\ge0}\frac{v_r(\eps,c)t^r}{r!}.
\]
The exponential factor removes the common Poisson normalization.
The next lemma explains why this configuration model generates the fixed-size quantities $\cZ_n$. A priori, the upper bounds in the lemma below could be infinite, in which case the inequalities would hold trivially. However, Lemmas~\ref{lem:series} and~\ref{lem:upper} show that these bounds are finite.
\begin{lemma}\label{lem:tile-product}
      For $t\ge0$, with $\cZ_0=1$,
      \begin{equation}\label{eq:poissonization}
            \mathcal Z(t)=\sum_{n\ge0}\frac{\cZ_n(\sigma,c)t^n}{n!}
            \le\prod_{\Delta'\in\cT}\mathcal V(tw_{\eps\Delta'}).
      \end{equation}
      Consequently, for $n\ge1$ and $t>0$,
      \begin{equation}\label{eq:occupancy}
            \cZ_n(\sigma,c)\le \frac{n!}{t^{n}}\prod_{\Delta'\in\cT}\cV(tw_{\eps\Delta'}).
      \end{equation}
\end{lemma}
\begin{proof}
      For $t>0$, conditional on having $n$ points, the process consists of $n$ independent points with density equal to the intensity divided by $t\Lambda$. By~\eqref{eq:Z-integral}, their antichain probability is $\cZ_n/\Lambda^n$. Multiplying by $e^{-t\Lambda}(t\Lambda)^n/n!$ and summing over $n$ proves the equality in~\eqref{eq:poissonization}; at $t=0$, both sides equal $1$.

      For a tile $B=\eps\Delta'$, let $\Pi_t|_B$ retain the points whose floors lie in $B$, and put $\Lambda_B=(1-2c)^{-1/2}\int_Bw$. The same expansion and Lemma~\ref{lem:tile-height} give
      \[
            e^{t\Lambda_B}\P(\Pi_t|_B\text{ is an antichain})
            =\sum_{r\ge0}\frac{t^r}{r!(2\pi)^{r/2}}
            \int_{\cA_r(B)}e^{-(\frac12-c)\|\fh\|^2}
            \prod_{k=1}^rw(\vf^{(k)})\dd\ff\dd\fh
            \le\mathcal V(tw_B).
      \]
      The term for $r=0$ is $1$. A global antichain restricts to an antichain in every tile. The restrictions of the Poisson process to distinct tiles are independent (since tile boundaries have zero intensity), and $\sum_B\Lambda_B=\Lambda$, so multiplication gives the product bound.
      Finally, every coefficient is nonnegative, so the $n$th term in the series is at most its sum, giving~\eqref{eq:occupancy}.
\end{proof}

It remains to bound the growth of the one-tile series, then sum its logarithm over the tiles.

\begin{lemma}\label{lem:series}
      For fixed $\eps>0$ and $c<1/2$, the exponential generating function $\cV(t)=\sum_{r\ge0}\frac{v_r(\eps,c)t^r}{r!}$ is finite for $t\ge0$, and
      \begin{equation}\label{eq:series}
            \limsup_{t\to\infty}\left(t^{-\frac{d-1}{d}} \log \cV(t)\right) \le\frac d{d-1}\exp\left(\frac{d-1}{d}b_\eps-\frac1d\right).
      \end{equation}
\end{lemma}
\begin{proof}
      Fix $\delta>0$. Equation~\eqref{eq:coefficient-rate} and Stirling's formula give a constant $C_\delta\ge1$ such that, for all integers $r\ge0$,
      \[
            \frac{v_r(\eps,c)}{r!} = (r!)^{-\frac d{d-1}} \left((r!)^{\frac{1}{d-1}}v_r(\eps,c)\right) \le  (r!)^{-\frac d{d-1}} C_\delta\exp\Bigl\{(b_\eps+\delta-\frac1{d-1})r\Bigr\}.
      \]
      For $t>0$, set
      \[
            s=\exp\left\{\frac{d-1}{d}(b_\eps+\delta)-\frac1d\right\} t^{\frac{d-1}{d}}.
      \]
      Then we have
      \[
            \cV(t) \le C_\delta\sum_{r\ge0}\left(\frac{s^r}{r!}\right)^{\frac d{d-1}}
            \le C_\delta\sum_{r\ge0}\Biggl(\frac{s^r}{r!}\Biggl(\sum_{j\ge0}\frac{s^j}{j!}\Biggr)^{\frac 1{d-1}}\Biggr)
            = C_\delta\Biggl(\sum_{r\ge0}\frac{s^r}{r!}\Biggr)^{\frac d{d-1}}
            =C_\delta e^{\frac d{d-1}s}.
      \]
      This proves convergence of the series. It follows that
      \[\limsup_{t\to\infty}\left(t^{-\frac{d-1}{d}} \log \cV(t)\right) \le\frac d{d-1}\exp\left(\frac{d-1}{d}(b_\eps+\delta)-\frac1d\right).\]
      Since this holds for any $\delta>0$, we obtain~\eqref{eq:series}.
\end{proof}

Finally we prove the upper bound for~\eqref{eq:integral-target}:
\begin{lemma}\label{lem:upper}
      For every bounded continuous $\sigma:E\to\R$ and $c<1/2$, with $w=\psi_Ee^\sigma$,
      \begin{equation}\label{eq:upper-bound}
            \limsup_{n\to\infty}\left(\frac1n\log\cZ_n(\sigma,c)+\frac{\log n}{d-1}\right)
            \le\gamma_d-\frac12\log(2\pi)+\frac d{d-1}\log\int_Ew(\vu)^{\frac{d-1}{d}}\dd\vu.
      \end{equation}
\end{lemma}
\begin{proof}
      For $0<\eps\le 1$, define the tilewise upper approximation $w_\eps:E\to(0,\infty)$ by
      \[
            w_\eps(\vu)=w_{\eps\Delta'}:=\sup\nolimits_{\eps\Delta'}w
            \qquad\text{for }\vu\in\eps\Delta',
      \]
      assigning each boundary point to an adjacent tile.
      We will show that \eqref{eq:upper-bound} holds with the integral replaced by
      \[I_\eps:=\int_E w_\eps(\vu)^{\frac{d-1}{d}}\dd\vu =
            |\eps\Delta|S_\eps,\qquad \text{where }
            S_\eps:=\sum_{\Delta'\in\cT}w_{\eps \Delta'}^{\frac{d-1}{d}}.
      \]
      We first show that $S_\eps<\infty$ and that $I_\eps\to\int_Ew(\vu)^{\frac{d-1}{d}}\dd\vu$ as $\eps\downarrow0$.

      Put $D=\diam\Delta$. If $0<\eps\le1$ and $\vu,\vv$ lie in the same tile $\eps\Delta'$, then $\|\vu-\vv\|\leq \eps D$, so
      \[\eps^2 D^2+\|\vu\|^2\ge \frac12(\eps D+\|\vu\|)^2\ge\frac{\|\vv\|^2}{2}.\] Since $\sigma$ is bounded, there are constants $C_1,C_2>0$, independent of $\eps$ and the tile, such that for every $0<\eps\le1$ and $\vv\in E$, with $\eps\Delta'$ denoting its assigned tile,
      \[
            w_\eps(\vv)=\sup_{\vu\in \eps\Delta'}(2\pi)^{-\frac{d-1}{2}}e^{\sigma(\vu)-\frac{\|\vu\|^{2}}{2}}\leq C_{1}\sup_{\vu\in \eps\Delta'}e^{-\frac{\|\vu\|^{2}}{2}}\leq C_{1} e^{-\frac{\|\vv\|^{2}}{4}+\frac{\eps^2 D^2}{2}}\leq C_2 e^{-\frac{\|\vv\|^{2}}{4}}.
      \]
      Raising this bound to the power $\frac{d-1}{d}$ gives the common integrable upper bound
      \begin{equation}\label{eq:tilewise-majorant}
            0<w_\eps(\vv)^{\frac{d-1}{d}}\le C_{2}^{\frac{d-1}{d}} e^{-\frac{(d-1)\|\vv\|^2}{4d}},\qquad \vv\in E.
      \end{equation}
      Hence $I_\eps<\infty$ and $S_\eps=I_\eps/|\eps\Delta|<\infty$. For each fixed $\vv\in E$, every point of the tile assigned to $\vv$ is within distance $\eps D$ of $\vv$. Continuity of $w$ therefore gives $w_\eps(\vv)\to w(\vv)$. From~\eqref{eq:tilewise-majorant}, the dominated convergence theorem yields, as $\eps\downarrow0$,
      \begin{equation}\label{eq:tilewise-integral-limit}
            |\eps\Delta|S_\eps=I_\eps\longrightarrow
            \int_Ew(\vu)^{\frac{d-1}{d}}\dd\vu\in(0,\infty).
      \end{equation}
      
      Now fix $\eps\in(0,1]$. By Lemma~\ref{lem:series}, together with $\log\cV(s)=O(s)$ as $s\downarrow0$, there is a constant $C=C(d,\eps,c)$ such that
      \[
            0\le t^{-\frac{d-1}{d}}\log\cV(tw_{\eps\Delta'})\le Cw_{\eps\Delta'}^{\frac{d-1}{d}},\qquad t>0,
      \]
      uniformly over the tiles. Since $S_\eps$ is finite, the right-hand side is summable, so the product in~\eqref{eq:occupancy} is finite. We may therefore sum the one-tile limsup bounds in~\eqref{eq:series}: retain finitely many tiles, then make the remaining sum uniformly small using the displayed summable bound. This gives
      \[
            \limsup_{t\to\infty}\biggl(t^{-\frac{d-1}{d}}\sum_{\Delta'\in\cT}\log\cV(tw_{\eps\Delta'})\biggr)
            \le K_\eps:=\frac d{d-1}e^{\frac{d-1}{d}b_\eps-\frac1d}S_\eps.
      \]
      Consequently, for fixed $\eps$,
      \[
            \log\cZ_n(\sigma,c) \le\log(n!)-n\log t+K_\eps t^{\frac{d-1}{d}} +o\bigl(t^{\frac{d-1}{d}}\bigr),\qquad t\to\infty.
      \]
      The remainder is independent of $n$. The leading expression $-n\log t+K_\eps t^{\frac{d-1}{d}}$ is minimized at
      \[
            t_n=t_n(\eps) :=\left(\frac{dn}{(d-1)K_\eps}\right)^{\frac d{d-1}},
      \]
      with $w,c$ still fixed. Since $t_n\to\infty$ and $t_n^{\frac{d-1}{d}}$ is proportional to $n$, the remainder is $o(n)$ at $t=t_n$. Stirling's formula now gives
      \[
            \begin{aligned}
                  \limsup_{n\to\infty} \left(\frac1n\log\cZ_n(\sigma,c)+\frac{\log n}{d-1}\right) & \le\frac d{d-1}\log\left(\frac{d-1}{d}K_\eps\right)+\frac1{d-1}                                              \\
                                                                                                  & =b_\eps+\frac d{d-1}\log S_\eps =\gamma_d-\frac12\log(2\pi) +\frac d{d-1}\log\bigl(|\eps\Delta|S_\eps\bigr).
            \end{aligned}
      \]
      Finally, let $\eps\downarrow0$ and apply~\eqref{eq:tilewise-integral-limit} to obtain~\eqref{eq:upper-bound}.
\end{proof}

\section{The lower bound for Proposition~\ref{prop:volume}}\label{sec:lower-bound}

The purpose of this section is to prove the lower bound for~\eqref{eq:integral-target}, namely~\eqref{eq:lower-conclusion}.
Lemma~\ref{lem:matching} constructs measurable sets $\cB_r\subset\cA_r(\Delta)$ such that independent copies of antichains in $\cB_r$ in different tiles $\Delta'$ always combine into an antichain. The volume of $\cB_r$ retains the logarithmic rate per point of $V_r(\Delta)$ as $r\to\infty$ and the contraction parameters tend to zero. We fill finitely many separated cubes with tiles, choosing the scale in each cube according to its floor weight, and then approximate the integral over $E$.

\subsection{Matching boundary values}

An antichain has a single height at each floor position. To combine independent choices in different tiles, we will select $\cB_r\subset\cA_r(\Delta)$ so that every configuration admits a height extension with the same boundary values on $\Delta$. The following proposition provides the extension used in this construction.

\begin{proposition}\label{prop:extension}
      Let $B\subset E$ be nonempty and bounded, and let $h:B\to\R$ be bounded. Suppose
      \begin{equation}\label{eq:height-data}
            h(\vu)-h(\vv)\le\Ldir{\vu-\vv} \qquad(\vu,\vv\in B).
      \end{equation}
      Then
      \[
            \widehat h(\vu):=\inf_{\vw\in B} \bigl(h(\vw)+\Ldir{\vu-\vw}\bigr),\qquad \vu\in E,
      \]
      is finite, agrees with $h$ on $B$, and satisfies
      \[
            \widehat h(\vu)-\widehat h(\vv)\le\Ldir{\vu-\vv} \qquad(\vu,\vv\in E).
      \]
      In particular, $\widehat h$ is $\sqrt{d-1}$-Lipschitz.
\end{proposition}
\begin{proof}
      Since $\mathcal L\ge0$, we have $\widehat h(\vu)\geq \inf_B h$, while any fixed $\vw\in B$ gives a finite upper bound on $\widehat h(\vu)$. Hence $\widehat h(\vu)$ is finite. If $\vu\in B$, then condition~\eqref{eq:height-data} implies $h(\vw)+\Ldir{\vu-\vw}\geq h(\vu)$ for each $\vw\in B$, while taking $\vw=\vu$ shows that $\widehat h(\vu)=h(\vu)$. For arbitrary $\vu,\vv\in E$, subadditivity of $\mathcal L$ gives
      \[
            h(\vw)+\Ldir{\vu-\vw} \le h(\vw)+\Ldir{\vv-\vw}+\Ldir{\vu-\vv} \qquad(\vw\in B).
      \]
      Taking the infimum over $\vw\in B$ proves the required bound. Finally, exchanging $\vu$ and $\vv$ gives $\widehat h(\vv)-\widehat h(\vu)\le\Ldir{\vv-\vu}$. Combining these two bounds with~\eqref{eq:height-bound} proves the stated Lipschitz constant.
\end{proof}

To construct $\cB_r$ and prove the following lemma, we first contract configurations toward an interior floor point to separate them from the boundary. We then contract the heights to leave enough room to replace nearby boundary functions by a common representative.

\begin{lemma}\label{lem:matching}
      Fix $\eta,\lambda\in (0,\frac12)$. There are $C=C(d)\ge 1$ and an integer $M=M(d,\eta,\lambda)\ge1$ such that, for every $r\in\N$, there exist a measurable set $\cB_r\subset\cA_r(\Delta)$ and a function $h_\partial:\partial\Delta\to\R$ with the following properties:
      \begin{enumerate}[label=\textup{(\roman*)},ref=\textup{(\roman*)},leftmargin=*]
            \item \label{item:fit-bounds}
                  Every $(\ff,\fh)\in\cB_r$ satisfies, for $1\le k\le r$,
                  \[
                        |h^{(k)}|\le C, \qquad \dist(\vf^{(k)},\partial\Delta)\ge\frac{\eta}{C}.
                  \]
            \item \label{item:fit-extension}
                  Every $(\ff,\fh)\in\cB_r$ admits an extension $\widehat h:\Delta\to\R$ such that
                  \[
                        \widehat h(\vf^{(k)})=h^{(k)},\qquad \widehat h|_{\partial\Delta}=h_\partial,\qquad \widehat h(\vu)-\widehat h(\vv)\le\Ldir{\vu-\vv} \quad(\vu,\vv\in\Delta).
                  \]
                  The extension may depend on the configuration, but its boundary values do not. The common boundary function may depend on $r$.
            \item \label{item:fit-volume}
                  The volume of $\cB_r$ satisfies
                  \begin{equation}\label{eq:matching}
                        \Vol_{dr}(\cB_r)\ge\frac{2}M(1-\eta)^{dr}(1-\lambda)^r V_r(\Delta).
                  \end{equation}
      \end{enumerate}
\end{lemma}

\begin{proof}
      Choose $\vf^\circ\in\Delta$ and $\rho_0>0$ such that the open ball of radius $\rho_0$ about $\vf^\circ$ is contained in $\Delta$.
      Take $(\ff,\mathbf a)\in\cA_r^0(\Delta)$ and $s\in[-1,1]$, and define $\fh$ by $h^{(k)}=a^{(k)}+s$, so $h^{(1)}=s$. The inverse map is then given by $s=h^{(1)}$, $a^{(k)}=h^{(k)}-h^{(1)}$. This map has Jacobian one in the free coordinates, so the starting family $\{(\ff,\fh):(\ff,\mathbf a)\in\cA_r^0(\Delta),s\in[-1,1]\}$ has full $dr$-dimensional volume $2V_r(\Delta)$. Apply
      \[
            \widetilde \vf^{(k)}=\eta\vf^\circ+(1-\eta)\vf^{(k)}, \qquad \widetilde h^{(k)}=(1-\eta)h^{(k)}.
      \]
      Let $(\widetilde \ff,\widetilde \fh)$ denote the image of $(\ff,\fh)$ under this transformation. By homogeneity of $\mathcal L$, $(\widetilde \ff,\widetilde \fh)$ is an antichain. Let $\cA_{r}^{\eta}(\Delta)$ denote the associated set of antichains
      \[\cA_{r}^{\eta}(\Delta)=\left\{(\widetilde \ff,\widetilde \fh):(\ff,\mathbf a)\in\cA_r^0(\Delta),\ s\in[-1,1],\ h^{(k)}=a^{(k)}+s\right\}.\]
      The Jacobian of this transformation is $(1-\eta)^{(d-1)r}(1-\eta)^r=(1-\eta)^{dr}$, so the volume $\Vol_{dr}(\cA_{r}^{\eta}(\Delta))=2(1-\eta)^{dr}V_r(\Delta)$. Moreover, each new floor point $\widetilde\vf^{(k)}$ has distance at least $\eta\rho_0$ from the boundary $\partial\Delta$. Indeed, for $\|\vu\|<\rho_0$, $\widetilde\vf^{(k)}+\eta\vu=(1-\eta)\vf^{(k)}+\eta(\vf^\circ+\vu)\in\Delta$, by the convexity of $\Delta$ since both $\vf^{(k)}\in\Delta$ and $\vf^\circ+\vu\in\Delta$.

      The bound $|h^{(1)}|=|s|\le1$ and Lemma~\ref{cor:height} give
      \[
            |\widetilde h^{(k)}|\le |h^{(k)}|\le |h^{(1)}|+|h^{(k)}-h^{(1)}|
            \le1+\sqrt{d-1}\diam(\Delta).
      \]
      The extension in Proposition~\ref{prop:extension} assigns to every $(\widetilde{\ff},\widetilde{\fh})\in\cA_r^\eta(\Delta)$ the function
      \begin{equation}\label{eq:def-ext}
            \widehat h_{\tilde\ff,\tilde\fh}(\vu)=\min_{1\le k\le r} \bigl(\widetilde{h}^{(k)}+\Ldir{\vu-\widetilde{\vf}^{(k)}}\bigr), \qquad \vu\in\Delta.
      \end{equation}
      These extensions are uniformly Lipschitz and bounded in absolute value by $1+2\sqrt{d-1}\diam(\Delta)$, independently of $r$. Fix $C=C(d)$ at least this large and at least $1/\rho_0$. Note then that the position and height bounds in part~\ref{item:fit-bounds} hold for any element of $\cA_{r}^{\eta}(\Delta)$.

      The boundary restrictions form a uniformly bounded, equicontinuous family of functions on $\partial\Delta$:
      \[
            \widehat{\cH}=\bigcup_{r=1}^{\infty}\left\{\widehat h_{\tilde\ff,\tilde\fh}\big|_{\partial\Delta}:(\widetilde{\ff},\widetilde{\fh})\in\cA_r^\eta(\Delta)\right\}.
      \]
      By the Arzel\`a--Ascoli theorem, this family is relatively compact in the uniform norm. Hence, for every $\tau>0$, there are finitely many representatives $h_\partial^{(1)},\ldots,h_\partial^{(M)}\in \widehat{\cH}$ such that every element of $\widehat{\cH}$ lies at distance less than $\tau$ from one of these representatives. Choosing $\tau=\frac{\lambda\eta\rho_0}{2\sqrt{d-1}}$, the value $M=M(d,\eta,\lambda)$ is independent of $r$ since the family $\widehat{\cH}$ does not depend on $r$. Each representative $h_\partial^{(j)}$ satisfies the directional height inequality~\eqref{eq:height-data}, because it is the restriction of an admissible extension.

      For fixed $r$, we partition the space $\cA_{r}^{\eta}(\Delta)$ by the first index $j$ for which $\|\widehat h_{\tilde\ff,\tilde\fh}|_{\partial\Delta}-h_\partial^{(j)}\|_\infty<\tau$. The boundary restriction depends continuously on $(\widetilde\ff,\widetilde\fh)$ in the uniform norm: each term in~\eqref{eq:def-ext} does so because $\mathcal L$ is Lipschitz, and taking a finite minimum preserves this continuity. Hence,
      these classes are measurable, and one of them has a volume of at least $\frac1M \Vol_{dr}(\cA_{r}^{\eta}(\Delta))$.

      Fix such a class $\cA_{r}^{(j)}(\Delta)$ and its associated representative $h_\partial^{(j)}$. Define $h_\partial=(1-\lambda)h_\partial^{(j)}$ on $\partial\Delta$ and define $\cB_r$ to be the class obtained from $\cA_{r}^{(j)}(\Delta)$ by multiplying all heights by $1-\lambda$. The height contraction has Jacobian $(1-\lambda)^r$, so
      \[\Vol_{dr}(\cB_r)=(1-\lambda)^r\Vol_{dr}(\cA_{r}^{(j)}(\Delta))\geq \frac1M(1-\lambda)^r\Vol_{dr}(\cA_{r}^{\eta}(\Delta))=\frac2M(1-\eta)^{dr}(1-\lambda)^rV_r(\Delta),\]
      which proves~\ref{item:fit-volume}. The bounds in~\ref{item:fit-bounds} hold since the same bounds hold for every element of $\cA_{r}^{\eta}(\Delta)$. It remains to show that~\ref{item:fit-extension} holds using our definitions of $\cB_r$ and $h_\partial$. Recall that the floors and heights of an antichain in $\cB_r$ can be written as $\widetilde\vf^{(k)},(1-\lambda)\widetilde h^{(k)}$. Put $\vu=\widetilde\vf^{(k)}$. Then for each $k$, each $\vv\in\partial\Delta$, and either choice of sign,
      \[
            \pm(1-\lambda)\bigl(\widetilde h^{(k)}-h_\partial^{(j)}(\vv)\bigr) \le(1-\lambda)\Ldir{\pm(\vu-\vv)}+(1-\lambda)\tau \le\Ldir{\pm(\vu-\vv)}.
      \]
      The last inequalities use $\Ldir{\pm(\vu-\vv)}\ge\frac{1}{\sqrt{d-1}}\|\vu-\vv\|\ge\frac{\eta\rho_0}{\sqrt{d-1}}$. In other words, the height contraction leaves slack $\lambda\Ldir{\pm(\vu-\vv)}$, which absorbs the term $(1-\lambda)\tau\le\lambda\eta\rho_0/(2\sqrt{d-1})$. The interior--interior and boundary--boundary inequalities also persist since $\mathcal L\ge0$. Hence we can apply Proposition~\ref{prop:extension} to the union of the finite floor set and $\partial\Delta$ to define an extension $\widehat{h}$ satisfying $\widehat{h}(\widetilde\vf^{(k)})=(1-\lambda)\widetilde h^{(k)}$ for each $k$ and $\widehat{h}|_{\partial\Delta}=h_{\partial}$, proving~\ref{item:fit-extension}. Note that $\cB_r$ is measurable because the height contraction is an invertible affine map and we already described why $\cA_{r}^{(j)}(\Delta)$ is measurable.
\end{proof}

\begin{lemma}\label{lem:replication}
      Fix $\eta,\lambda\in(0,\frac12)$ and $r\in\N$, and let $C$, $\cB_r$, and $h_\partial$ be as in Lemma~\ref{lem:matching}. Choose $q\ge1$ distinct tiles at a common scale $\eps>0$, and put $N=rq$. In each tile $\eps\Delta'=\eps(\omega\Delta)$ where $\omega$ is a composition of side reflections, copy a configuration from $\cB_r$ by
      \[
            (\vf^{(k)},h^{(k)})\longmapsto(\eps\,\omega(\vf^{(k)}),\eps h^{(k)}).
      \]
      Independent choices in the tiles give an admissible $N$-point configuration, with every height bounded in absolute value by $C\eps$. Assigning an unordered $r$-element label set to each tile, using increasing label order within that tile, gives a measurable family of volume
      \begin{equation}\label{eq:replication}
            \frac{N!}{(r!)^q}\bigl(\eps^{dr}\Vol_{dr}(\cB_r)\bigr)^q.
      \end{equation}
\end{lemma}
\begin{proof}
      Transport a height extension from Lemma~\ref{lem:matching} to each chosen tile, and use any such extension on unused tiles. We will show that this defines a global function $h:E\to \R$ satisfying~\eqref{eq:height-data}.
      By Proposition~\ref{prop:tiling}, the transported copies of the common boundary values agree wherever tiles meet, so $h$ is well defined. Since scaling and coordinate permutations preserve the directional inequalities, \eqref{eq:height-data} holds for any two points in the same tile. For any $\vu,\vv\in E$, local finiteness and convexity of the tiles give a partition $\vu=\vu_0,\vu_1,\ldots,\vu_m=\vv$ of the segment $[\vu,\vv]$ into subsegments lying in single tiles. Hence
      \[
            h(\vu)-h(\vv)=\sum_{k=0}^{m-1}\bigl(h(\vu_k)-h(\vu_{k+1})\bigr)\le\sum_{k=0}^{m-1}\Ldir{\vu_k-\vu_{k+1}}=\Ldir{\vu-\vv}.
      \]
      The last equality holds because each displacement is a nonnegative multiple of $\vu-\vv$ and $\mathcal L$ is positively homogeneous.
      Hence the $N$ points satisfy~\eqref{eq:constraints}, so they form an antichain.
      From~\ref{item:fit-bounds} in Lemma~\ref{lem:matching}, the sampled heights are bounded by $C\eps$.

      There are $N!/(r!)^q$ ways to assign the $N$ labels to the $q$ tiles, with $r$ labels assigned to each tile. The floor scaling by $\eps$ has Jacobian $\eps^{(d-1)r}$ and the height scaling has Jacobian $\eps^r$, so for each such assignment, each tile contributes volume $\eps^{dr}\Vol_{dr}(\cB_r)$. Each label assignment therefore contributes volume $(\eps^{dr}\Vol_{dr}(\cB_r))^{q}$. The resulting configuration sets are disjoint because all floor points lie in tile interiors, so the total volume is~\eqref{eq:replication}.
\end{proof}

\subsection{Weighted lower bound}
Fix a bounded continuous function $\sigma:E\to\R$ and $c<1/2$.
We use finitely many cubes $\square_j$, placing $r$ points in each complete simplex $\eps_j\Delta'$ inside a cube. We choose the tile scales $\eps_{j}$ using the value of $w=\psi_Ee^\sigma$ in each cube $\square_{j}$.

\begin{lemma}\label{lem:cube-lower}
      Let $\square_1,\ldots,\square_J\subset E$ be closed cubes of positive volume with disjoint interiors. Set
      \[
            w_j=\inf_{\square_j}w>0,\qquad S=\sum_{j=1}^J|\square_j|w_j^{\frac{d-1}{d}}.
      \]
      Then
      \begin{equation}\label{eq:lower-rate}
            \liminf_{n\to\infty}\left(\frac1n\log\cZ_n(\sigma,c)+\frac{\log n}{d-1}\right)
            \ge\gamma_d-\frac12\log(2\pi)+\frac d{d-1}\log S.
      \end{equation}
\end{lemma}
\begin{proof}
      First suppose that the cubes are positively separated. Fix $\eta,\lambda\in(0,\frac12)$ and $r\in\N$, and choose $C$, $M=M(d,\eta,\lambda)$, and $\cB_r$ as in Lemma~\ref{lem:matching}. For $\eps>0$, set $\eps_j=\eps w_j^{-1/d}$, so $w_j\eps_j^d=\eps^d$ for every $j$. In each complete tile of scale $\eps_j$ contained in $\square_j$, place a scaled copy of an independently chosen configuration from $\cB_r$. Let $q_j=q(\eps_j,\square_j)$ be the number of these tiles, and put
      \[
            N_j=rq_j,\qquad N=\sum_jN_j.
      \]
      Let $\cA_N^\eps(E)$ denote the family of all labeled configurations obtained in this way. We check that $\cA_N^\eps(E)\subset\cA_N(E)$ for all sufficiently small $\eps$. Lemma~\ref{lem:replication} gives admissibility within each cube. If $J>1$, set $a=\min_{i\ne j}\dist(\square_i,\square_j)>0$. For points in different cubes, by Lemma~\ref{lem:matching},
      \[
            |h-h'|\le2C\max_j\eps_j\le a/\sqrt{d-1}
            \le\Ldir{\pm(\vf-\vf')}
      \]
      for small enough $\eps$. For $J=1$ only the within-cube check is needed. Hence $\cA_N^\eps(E)\subset\cA_N(E)$ for all sufficiently small $\eps$.

      Proposition~\ref{prop:tiling} gives
      \begin{equation}\label{eq:qj_asymptotics}
            q_j=\frac{|\square_j|}{\eps_j^{d-1}|\Delta|}+O(\eps_j^{2-d})
            \qquad(\eps_j\downarrow0).
      \end{equation}
      Since $N_j=rq_j$ and $N=\sum_jN_j$, substituting $\eps_j=\eps w_j^{-1/d}$ gives
      \begin{equation}\label{eq:N_Nj_asymptotics}
            N_j\eps^{d-1}=\frac{r|\square_j|w_j^{\frac{d-1}{d}}}{|\Delta|}+O(\eps),
            \qquad
            N\eps^{d-1}=\frac{rS}{|\Delta|}+O(\eps).
      \end{equation}

      Since the integrand in~\eqref{eq:Z-integral} is nonnegative, restricting it to $\cA_N^\eps(E)$ gives a lower bound for $\cZ_N(\sigma,c)$. Every height has absolute value at most $C\max_j\eps_j$, so the height factor is at least $e^{-N\kappa_\eps}$, where $\kappa_\eps=(\frac12-c)C^2(\max_j\eps_j)^2\to0$. There are $N!/\prod_jN_j!$ ways to allocate the labels to the cubes. For each allocation, the choices in different cubes can be made independently, so their volumes multiply. Lemma~\ref{lem:replication} gives the configuration volume in cube $\square_j$ as
      \[
            \frac{N_j!}{(r!)^{q_j}}\bigl(\eps_j^{dr}\Vol_{dr}(\cB_r)\bigr)^{q_j}.
      \]
      Since the floor weight is at least $w_j$ in this cube,
      \begin{equation}\label{eq:lower-volume}
            \begin{aligned}
                  \cZ_N(\sigma,c)
                   & \ge e^{-N\kappa_\eps}(2\pi)^{-N/2}N!
                  \prod_j\frac{\bigl(w_j^r\eps_j^{dr}\Vol_{dr}(\cB_r)\bigr)^{q_j}}{(r!)^{q_j}} \\
                   & =N!e^{-N\kappa_\eps}(2\pi)^{-N/2}
                  \left(\frac{\Vol_{dr}(\cB_r)}{r!}\right)^{N/r}\eps^{dN}.
            \end{aligned}
      \end{equation}
      The equality uses $w_j\eps_j^d=\eps^d$ and $\sum_jq_j=N/r$. Here and below, $o_\eps(1)$ tends to zero as $\eps\downarrow0$, with $r,\eta,\lambda$ and the cubes fixed. Taking logarithms in~\eqref{eq:lower-volume} and applying Stirling's formula to $N!$ gives
      \[
            \frac1N\log\cZ_N(\sigma,c)+\frac{\log N}{d-1}
            \ge{} \frac1r\log\frac{\Vol_{dr}(\cB_r)}{r!}
            -1-\frac12\log(2\pi)
            +\frac d{d-1}\log(N\eps^{d-1})-o_\eps(1),
      \]
      where the identity $\log N+\log N/(d-1)+d\log\eps=\frac d{d-1}\log(N\eps^{d-1})$ collects the size and scale terms.
      For the local factor, \eqref{eq:matching} gives
      \[
            \frac1r\log\frac{\Vol_{dr}(\cB_r)}{r!}
            \ge\frac1r\log V_r(\Delta)+d\log(1-\eta)+\log(1-\lambda)
            -\frac{\log(M/2)+\log(r!)}r.
      \]
      Set
      \[
            e_r:=\frac{\log(M/2)+\log(r!)-r\log r+r}{r}
            =\frac{\log(M/2)}r+O\!\left(\frac{\log(r+1)}r\right).
      \]
      Substituting~\eqref{eq:gamma} yields
      \[
            \frac1r\log\frac{\Vol_{dr}(\cB_r)}{r!}
            \ge{} \gamma_{d,r}-\frac d{d-1}\log r
            +\frac{d\log|\Delta|}{d-1}+1
            +d\log(1-\eta)+\log(1-\lambda)-e_r.
      \]
      The error $e_r$ is independent of $\eps$ and tends to zero as $r\to\infty$ with $\eta,\lambda$ fixed.

      To pass from the constructed sizes $N=N(\eps)$ to a prescribed $n$, keep $r,\eta,\lambda$ fixed and choose the common scale so that $N$ is slightly larger than $n$. Take
      \[
            1-\frac1{d-1}<\zeta<1,
            \qquad
            \eps=\left(\frac{rS}{|\Delta|(n+n^\zeta)}\right)^{\frac1{d-1}}.
      \]
      Since complete tiles occupy at most the volume of each cube, \eqref{eq:N_Nj_asymptotics} gives
      \[
            0\le n+n^\zeta-N
            =O(\eps^{2-d})
            =O\bigl(n^{1-\frac1{d-1}}\bigr).
      \]
      Since $\zeta>1-\frac1{d-1}$, we have $n\le N\le n+n^\zeta$ for all sufficiently large $n$; since $\zeta<1$, we also have $(N-n)\log n=o(n)$. To compare these sample sizes directly in~\eqref{eq:Z-integral}, observe that the first $n$ points of an $N$-point antichain still form an antichain. Drop all constraints involving the remaining $N-n$ points and integrate those points freely. Each contributes the one-point mass in~\eqref{eq:poisson-mass}, so
      \begin{equation*}
            \cZ_n(\sigma,c)\ge
            \left(\frac{\int_Ew(\vu)\dd\vu}{\sqrt{1-2c}}\right)^{-(N-n)}
            \cZ_N(\sigma,c).
      \end{equation*}
      The logarithm of this prefactor and $N\log N-n\log n$ are $o(n)$, while $N/n\to1$, so this comparison transfers the normalized lower bound from $N$ to $n$. Taking $n\to\infty$ at fixed $r,\eta,\lambda$ removes the $o_\eps(1)$ terms. Substituting the local-factor bound and the limit of $N\eps^{d-1}$ gives
      \[
            \begin{aligned}
                  \liminf_{n\to\infty}
                  \left(\frac1n\log\cZ_n(\sigma,c)+\frac{\log n}{d-1}\right)
                  \ge{} & \gamma_{d,r}-\frac12\log(2\pi)
                  +\frac d{d-1}\log S                         \\
                        & +d\log(1-\eta)+\log(1-\lambda)-e_r.
            \end{aligned}
      \]
      The left-hand side is independent of $r$. Let $r\to\infty$ along a subsequence for which $\gamma_{d,r}\to\gamma_d$, and then let $\eta,\lambda\downarrow0$. This proves~\eqref{eq:lower-rate} for positively separated cubes.

      For cubes with only disjoint interiors, apply the bound to their concentric contractions by a factor $s\in(0,1)$. These are compact subsets of the respective disjoint interiors. The contracted cubes are positively separated and their lower sum is at least $s^{d-1}S$, since their weight infima can only increase. Letting $s\uparrow1$ proves the same bound.
\end{proof}

\begin{lemma}\label{lem:lower}
      For every bounded continuous $\sigma:E\to\R$ and $c<1/2$, with $w=\psi_Ee^\sigma$,
      \begin{equation}\label{eq:lower-conclusion}
            \liminf_{n\to\infty}\left(\frac1n\log\cZ_n(\sigma,c)+\frac{\log n}{d-1}\right)
            \ge\gamma_d-\frac12\log(2\pi)+\frac d{d-1}\log\int_Ew(\vu)^{\frac{d-1}{d}}\dd\vu.
      \end{equation}
\end{lemma}
\begin{proof}
      Fix $\delta>0$. Since $\sigma$ is bounded, \eqref{eq:floor_weight_def} gives $w(\vu)^{\frac{d-1}{d}}\le C e^{-\frac{d-1}{2d}\|\vu\|^2}$ for some $C>0$, so the integral is finite. Choose a large cube with tail integral below $\delta/2$, then use uniform continuity of $w^{(d-1)/d}$ within the cube to obtain a lower Riemann sum for the integral with error below $\delta/2$. This yields finitely many cubes with pairwise disjoint interiors such that
      \[
            S:=\sum_{j=1}^J|\square_j|w_j^{\frac{d-1}{d}}\ge\int_Ew(\vu)^{\frac{d-1}{d}}\dd\vu-\delta,
      \]
      where $w_j=\inf_{\square_j}w>0$ for each $j$. Applying Lemma~\ref{lem:cube-lower} and letting $\delta\downarrow0$ yields~\eqref{eq:lower-conclusion}.
\end{proof}

\begin{proof}[Proof of Proposition~\ref{prop:volume}]
      The upper bound~\eqref{eq:upper-bound} and the lower bound~\eqref{eq:lower-conclusion} establish~\eqref{eq:integral-target}. This gives~\eqref{eq:weighted} upon substituting $w=\psi_Ee^\sigma$.
\end{proof}

\section{Applications and open questions}\label{sec:applications}
In this section, we discuss several applications of our results to limits and enumeration of pattern-avoiding permutations and intervals in the weak Bruhat order. We conclude with some related open questions.

\subsection{Permutons and enumeration}\label{sec:ranks}

Let $\mathfrak S_n$ denote the set of permutations of $[n]=\{1,\ldots,n\}$. Before conditioning on $\cA_n$, the points $\vX^{(1)},\ldots,\vX^{(n)}$ are independent and uniform on $[0,1]^d$. For each coordinate $i$, write $S_{i,1}<\cdots<S_{i,n}$ for the values $\vX_i^{(1)},\ldots,\vX_i^{(n)}$ in increasing order; coordinate ties have probability zero.
Order the sample by its first coordinate, and let $\pi_i(k)$ be the rank of the $i$th coordinate of the $k$th point among these $n$ values, so $\pi_1=\mathrm{id}$. The coordinate samples $\vX_i^{(j)}$ are independent. By symmetry, conditional on the sorted values $S_{i,1}<\cdots<S_{i,n}$ in each coordinate, all $n!$ orders of the sample labels are equally likely.
Thus the rank tuple is uniform over the $(n!)^{d-1}$ possible $d$-permutations $(\mathrm{id},\pi_2,\ldots,\pi_d)$, with $\pi_i\in\mathfrak S_n$, and independent of the full collection of sorted coordinate values $(S_{i,r})_{i,r}$.

The antichain event depends only on the ranks, so conditioning on $\cA_n$ makes the corresponding rank tuple uniform on the set $\cC_n$ of $d$-permutations avoiding the pattern $(12,12,\ldots,12)$; see Section~\ref{sec:background} for definitions. Recalling $a_n=\P(\cA_n)$, we obtain
\begin{equation}\label{eq:rank-count}
      |\cC_n|=(n!)^{d-1}a_n.
\end{equation}
As a consequence, Theorem~\ref{thm:counting} translates to the following enumerative result:
\begin{proposition}\label{prop:counting}\label{cor:counting}
      Set $\beta_d'=\exp\!\left(\beta_d-\frac1{d-1}\right)>0$. Then
      \[
            |\cC_n|=(n!)^{d-1-\frac1{d-1}}(\beta_d')^n e^{o(n)}.
      \]
      The same asymptotic formula holds for avoidance of any prescribed two-point $d$-dimensional permutation pattern. In particular, $\beta_2'=1$.
\end{proposition}
\begin{proof}
      Equation~\eqref{eq:rank-count}, Theorem~\ref{thm:counting}, and Stirling's formula give
      \[
            \log|\cC_n|-\left(d-1-\frac1{d-1}\right)\log(n!) =n\left(\beta_d-\frac1{d-1}\right)+o(n).
      \]
      Reversing suitable coordinate axes gives bijections between all two-point avoidance classes. For $d=2$, the only allowed second permutation is reversal, so $|\cC_n|=1$.
\end{proof}
The earlier estimates~\cite{CK,GP} determine the factorial scale up to exponential factors. The values of the growth constants for $d\ge3$ and the error term are not determined here; see Question~\ref{ques:constant}. Nonetheless, this answers the two-point case of~\cite[Question~9.3]{CK}. The more general question asks whether the same formula (with a pattern-dependent positive constant in place of $\beta_d'$) holds for avoidance of any $d$-dimensional pattern of length at least $2$.

We now describe the permuton limit of $(12,\ldots,12)$-avoiding $d$-permutations $(\pi_1,\ldots,\pi_d)$.
Write
\[
      u_{r,n}=\frac{r-\tfrac12}{n},\qquad \vR^{(k)}=(u_{\pi_1(k),n},\ldots,u_{\pi_d(k),n}),\qquad \widehat\mu_n^{\mathrm{rk}}=\frac1n\sum_{k=1}^n\delta_{\vR^{(k)}}.
\]
Following~\cite{BL}, define $\mu_n^{\mathrm{rk}}$ by spreading each point mass $\frac1n$ at $\vR^{(k)}$ of $\widehat\mu_n^{\mathrm{rk}}$ uniformly over the cube of side $\frac1n$ centered at that rank point $\vR^{(k)}$. In each coordinate, every rank interval $[\frac{r-1}{n},\frac{r}{n}]$ occurs once and receives mass $\frac1n$, so $\mu_n^{\mathrm{rk}}$ has uniform marginals and is a permuton. Recall that $\mathcal P_d$ denotes the probability measures on $[0,1]^d$ with the weak topology, and that $\mu^\ast=\law(\Phi(\vF^\ast))$ is the limiting measure in~\eqref{eq:mu_n_hat_definition}.

      \begin{theorem}[The antichain permuton]\label{thm:permuton}\label{prop:permuton}\label{cor:permuton}
            For every open neighborhood $\cU$ of $\mu^\ast$ in $\mathcal P_d$, there is $\alpha=\alpha(d,\cU)>0$ such that
            \[
                  \P_n\bigl(\widehat\mu_n^{\mathrm{rk}}\notin\cU \text{ or }\mu_n^{\mathrm{rk}}\notin\cU\bigr) \le e^{-\alpha n}
            \]
            for all sufficiently large $n$.
      \end{theorem}
      \begin{proof}
            Conditioned on $\cA_n$, the sorted values $S_{i,1},\ldots,S_{i,n}$ in each coordinate still have the law of $n$ independent uniform samples arranged in increasing order, by the rank independence. Relabel the original points in increasing order of their first coordinate as $\widetilde{\vX}^{(1)},\ldots,\widetilde{\vX}^{(n)}$. Then
            \[
                  \widetilde\vX^{(k)}=(S_{1,k},S_{2,\pi_2(k)},\ldots,S_{d,\pi_d(k)}).
            \]
            The reordered sample has the same empirical measure $\widehat\mu_n$, as defined in~\eqref{eq:mu_n_hat_definition}.

            For $t\in[0,1]$, put $N_i(t)=\#\{r:S_{i,r}\le t\}$. Conditioned on $\cA_n$, this count has distribution $\operatorname{Binomial}(n,t)$. Fix $s>0$. If $S_{i,r}>u_{r,n}+s$, then $N_i(u_{r,n}+s)\le r-1$, at least $ns$ below its mean $r-\frac12+ns$. Similarly, $S_{i,r}<u_{r,n}-s$ forces $N_i(u_{r,n}-s)\ge r$, at least $ns$ above its mean. Hoeffding's inequality~\cite[Theorem~1]{Hoeffding1963} therefore gives
            \[
                  \P_n\bigl(|S_{i,r}-u_{r,n}|>s\bigr)\le2e^{-2ns^2}.
            \]
            Pairing each reordered sample point $\widetilde{\vX}^{(k)}$ with its rank point $\vR^{(k)}$ gives
            \[
                  D_n:=\max_k\|\widetilde\vX^{(k)}-\vR^{(k)}\|
                  \le\sqrt d\max_{i,r}|S_{i,r}-u_{r,n}|,
            \]
            so a union bound yields
            \begin{equation}\label{eq:rank-displacement}
                  \P_n(D_n>s)\le2dn\exp\left(-\frac{2ns^2}{d}\right).
            \end{equation}

            We combine~\eqref{eq:rank-displacement} with Theorem~\ref{thm:empirical-limit} by comparing finitely many continuous test functions. Choose continuous functions $g_1,\ldots,g_m$ and $\eps>0$ such that
            \[
                  \max_j\left|\int g_j\dd\nu-\int g_j\dd\mu^\ast\right|<\eps
                  \quad\Longrightarrow\quad\nu\in\cU.
            \]
            By uniform continuity, choose $s>0$ so that moving a point by at most $s$ changes each $g_j$ by less than $\eps/3$. Then on $\{D_n\le s\}$, replacing the sampled points $\widetilde\vX^{(k)}$ by rank points $\vR^{(k)}$ changes each test average
            $\frac{1}{n}\sum_{k=1}^{n}g_{j}(\widetilde\vX^{(k)})$
            by less than $\eps/3$. Spreading over the rank boxes changes it by less than another $\eps/3$ for large $n$, since every point of a box lies within $\sqrt d/(2n)$ of its center. Thus both rank measures $\widehat\mu_n^{\mathrm{rk}}$ and $\mu_n^{\mathrm{rk}}$ lie in $\cU$ whenever $D_n\le s$ and $\left|\int g_j\dd\widehat\mu_n-\int g_j\dd\mu^\ast\right|<\eps/3$ for every $j$. Theorem~\ref{thm:empirical-limit} and~\eqref{eq:rank-displacement} bound the failure probability exponentially in $n$.
      \end{proof}

\begin{remark}\label{rem:rank-hausdorff}
      For the rank point cloud $A_n^{\mathrm{rk}}=\{\vR^{(1)},\ldots,\vR^{(n)}\}$, matching points and applying the triangle inequality gives
      \[
            d_{\rm H}(A_n^{\mathrm{rk}},\overline{\cS_d})
            \le D_n+d_{\rm H}(A_n,\overline{\cS_d}).
      \]
      For each $\delta>0$, \eqref{eq:rank-displacement} and Theorem~\ref{thm:shape} bound the probability that this distance exceeds $\delta$ by $e^{-\alpha n}$, for some $\alpha=\alpha(d,\delta)>0$ and all sufficiently large $n$.
\end{remark}

\begin{remark}\label{rem:pattern-concentration}
      Theorem~\ref{thm:permuton} also gives exponential concentration of every fixed pattern frequency around its probability under independent sampling from $\mu^\ast$. Pattern probabilities are continuous in the permuton, and for a pattern on $q$ points the permutation frequency differs from its permuton probability by at most $\binom{q}{2}/n$; see~\cite[Lemma~2.3 and Proposition~2.4]{BL}.
\end{remark}


\subsection{Intervals in the weak Bruhat order}

In this section, set $d=3$, so a rank tuple in $\cC_n$ is $(\mathrm{id},\pi_2,\pi_3)$. For $\pi\in\mathfrak S_n$, define the inversion set of $\pi$ by
\[
      \Inv(\pi)=\{(k,\ell):k<\ell,\ \pi(k)>\pi(\ell)\}.
\]
The {\em weak Bruhat order} $\le_{\rm W}$ on $\mathfrak S_n$ is defined by $\pi^-\le_{\rm W}\pi^+$ when $\Inv(\pi^-)\subseteq\Inv(\pi^+)$. Thus moving upward in this order adds inversions. The interval
\[
      [\pi^-,\pi^+]_{\rm W}=\{\pi\in\mathfrak S_n:\pi^-\le_{\rm W}\pi\le_{\rm W}\pi^+\}
\]
is determined by its lower and upper endpoints $\pi^-$ and $\pi^+$. With $w_0(r)=n+1-r$, the correspondence with simultaneous pattern avoidance takes the form below; see~\cite{GP} and~\cite[Proposition~5]{BM}:
\begin{equation}\label{eq:interval-map}
      (\mathrm{id},\pi_2,\pi_3)\longmapsto (\pi^-,\pi^+)=(w_0\circ\pi_3,\pi_2).
\end{equation}
Indeed,
\[
      \Inv(w_0\circ\pi_3)=\{(k,\ell)\in[n]^2:k<\ell\}\setminus\Inv(\pi_3).
\]
The antichain condition says precisely that every increasing pair of $\pi_3$, hence every inversion of $w_0\circ\pi_3$, is an inversion of $\pi_2$. The inverse map is $(\pi_2,\pi_3)=(\pi^+,w_0\circ\pi^-)$, so a uniform tuple in $\cC_n$ gives a uniform weak-order interval.

\begin{proposition}
      Let $I_n$ be the number of weak-order intervals in $\mathfrak S_n$, and let $\Pi_1,\Pi_2$ be independent uniform permutations in $\mathfrak S_n$. Then
      \[
            I_n=(n!)^{3/2}(\beta_3')^n e^{o(n)},\qquad \P(\Pi_1\le_{\rm W}\Pi_2)=(n!)^{-1/2}(\beta_3')^n e^{o(n)}.
      \]
\end{proposition}

\begin{proof}
      By~\eqref{eq:interval-map}, $I_n=|\cC_n|$, while the probability is $I_n/(n!)^2$. Apply Proposition~\ref{cor:counting} with $d=3$.
\end{proof}

This probability was first studied by Hammett and Pittel~\cite{HP}. The leading asymptotic $\exp((-\frac12+o(1))n\log n)$ is given in~\cite[Theorem~1.1]{BCHHLM} and also follows from the earlier coordinate-order and pattern-avoidance bounds~\cite{CK,GP}.
The permuton limit also describes the joint geometry of the interval endpoints.

\begin{proposition}\label{prop:intervals}
      Let $(\pi^-,\pi^+)$ be a uniformly chosen weak-order interval in $\mathfrak S_n$. If $(A,B)$ is a centered bivariate Gaussian with variances one and covariance $-1/2$, then
      \begin{equation}\label{eq:interval-law}
            \frac1n\sum_{k=1}^n \delta_{(u_{k,n},u_{\pi^-(k),n},u_{\pi^+(k),n})} \longrightarrow \law\bigl(\Phi(A),\Phi(A+B),\Phi(B)\bigr)
      \end{equation}
      in probability, with exponential concentration outside every fixed weak neighborhood.
\end{proposition}
\begin{proof}
      For $d=3$, $\vF^\ast$ defined in~\eqref{eq:Fstar_definition} is distributed as $(A,B,-A-B)$. The map $\vx\mapsto(\vx_1,1-\vx_3,\vx_2)$ implements~\eqref{eq:interval-map}, so Theorem~\ref{thm:permuton} gives~\eqref{eq:interval-law}.
\end{proof}

\begin{corollary}\label{cor:interval-inversions}
      For $n\ge2$, let $(\pi^-,\pi^+)$ be a uniformly chosen weak-order interval in $\mathfrak S_n$, and write $\inv(\pi):=|\Inv(\pi)|/\binom n2$ for the inversion density of $\pi\in\mathfrak S_n$. For every $\delta>0$, there is $\alpha=\alpha(\delta)>0$ such that
      \[
            \P\left(\max\left\{
            \Bigl|\inv(\pi^-)-\frac13\Bigr|, 
            \Bigl|\inv(\pi^+)-\frac23\Bigr|\right\}
            >\delta\right)\le e^{-\alpha n}
      \]
      for all sufficiently large $n$.
\end{corollary}
\begin{proof}
      By~\eqref{eq:interval-map}, inversions of $\pi^-$ are increasing pairs of $\pi_3$. Such a pair must be decreasing in $\pi_2$, since $(12,12,12)$ is forbidden by the antichain condition. Thus $\inv(\pi^-)$ is the frequency of $(12,21,12)$ in the antichain rank tuple $(\mathrm{id},\pi_2,\pi_3)$. An inversion of $\pi^+=\pi_2$ gives $(12,21,12)$ or $(12,21,21)$ according as $\pi_3$ increases or decreases along that pair, so $\inv(\pi^+)$ is the sum of their frequencies.

      Two independent samples from $\mu^\ast$ lie in the antichain $\cS_3$ and have no coordinate ties almost surely. Their pattern is therefore one of $(12,21,12)$, $(12,21,21)$, and $(12,12,21)$. These three patterns correspond to the three choices of the coordinate whose ordering differs from the other two. Since $\mu^\ast$ is invariant under coordinate permutations, each pattern has probability $1/3$. Remark~\ref{rem:pattern-concentration} therefore gives exponential concentration of the two inversion densities around $1/3$ and $2/3$, respectively.
\end{proof}



\subsection{Open questions}

Theorem~\ref{thm:counting} determines the logarithmic counting asymptotic up to an $o(n)$ error, but does not evaluate the constant $\beta_d$ for $d\ge3$. For comparison, when $d=2$ the identity $a_n=1/n!$ gives $\beta_2=1$ and remainder $-\frac12\log(2\pi n)+O(\frac1n)$ by Stirling's formula.

\begin{question}\label{ques:constant}
            What are the values of $\beta_d$, or equivalently $\gamma_d$, for $d\ge3$? Additionally, determine the order of the remainder
            \[
                  \log a_n+\frac{n}{d-1}\log n-\beta_dn.
            \]
            In particular, what is the weak-order interval constant $\beta_3'$, and is the remainder $O(\sqrt n)$ when $d=3$?
      \end{question}

Recall the total squared height $T_n=\sum_{k=1}^n(H^{(k)})^2$. The proof of Theorem~\ref{thm:empirical-limit} gives exponential concentration of $T_n/n$ near zero, while Theorem~\ref{thm:shape} controls outliers and holes. Neither result identifies the scale of the geometric deviations.

\begin{question}
      What are the typical orders of $T_n$ and $d_{\rm H}(A_n,\overline{\cS_d})$ under $\P_n$? Do these quantities have nondegenerate limiting distributions after suitable centering and rescaling?
\end{question}

For $g(\vx)=\vx_i$, the rank independence in Section~\ref{sec:ranks} leaves the empirical average's law unchanged by conditioning. Thus $\sqrt n(\int g\dd\widehat\mu_n-\frac12)$ converges in distribution to $\mathcal N(0,1/12)$. For general smooth test functions, the scale and limiting law of empirical fluctuations remain open. Recall that $\E_n=\E(\,\cdot\mid\cA_n)$.

\begin{question}
      For fixed $d\ge3$ and a smooth function $g:[0,1]^d\to\R$, determine the normalization and limiting distribution of
      \[
            \int g\dd\widehat\mu_n
            -\E_n\left[\int g\dd\widehat\mu_n\right].
      \]
      For which such $g$ is the limiting distribution Gaussian? 
\end{question}

Coordinate averages of $\widehat\mu_n^{\mathrm{rk}}$ are deterministic, unlike those of $\widehat\mu_n$, so the two measures need not have the same fluctuations.

For an integer $K\ge1$, let $\cC_{n,K}$ denote the set of $d$-permutations of $[n]$ avoiding the increasing pattern of size $K+1$; equivalently, their rank point set contains no chain of length $K+1$ under $\prec$. Let $\mu_{n,K}^{\mathrm{rk}}$ be the rank-box permuton of a uniform member of $\cC_{n,K}$, as in Section~\ref{sec:ranks}. Coloring each point by the length of a longest chain ending there gives a $K$-coloring, and each color class forms an antichain. This suggests the following extension of our results.
\begin{conjecture}\label{conj:bounded-chains} For every fixed $d\ge3$ and $K\ge2$, as $n\to\infty$,
      \begin{equation}\label{eq:bounded-chains-counting}
            |\cC_{n,K}|=(n!)^{d-1-\frac1{d-1}}
            \left(\beta_d' K^{\frac d{d-1}}\right)^n e^{o(n)},
      \end{equation}
      where $\beta_d'$ is the counting constant in Proposition~\ref{cor:counting}. Moreover, the corresponding permuton $\mu_{n,K}^{\mathrm{rk}}$ converges weakly in probability to $\mu^\ast$.
\end{conjecture}

This leads to a question regarding the universality of pattern-avoiding permutations.
\begin{question}\label{ques:universality} For which patterns $\tau$ of size at least $3$ is $\mu^\ast$ (or its reflection in some coordinate axes) the permuton limit of uniform $\tau$-avoiding $d$-permutations?
\end{question}
The same question can be asked for avoidance of several patterns, provided the classes are nonempty for every $n$. A general classification remains open even for $d=2$; see~\cite{hohmeier2025permuton}.

\section*{Acknowledgments}
Part of this work was carried out during visits hosted by Ewain Gwynne, with support from National Science Foundation grant DMS-2245832. We thank Richard Kenyon and Istv\'an Prause for sharing their insights into surface tension related to our results, and Jacopo Borga for sharing a draft survey on permutons.

\section*{AI disclosure statement}\label{sec:ai}
\addcontentsline{toc}{section}{AI disclosure statement}
This work was initiated when we posed the problem for $d=3$ to OpenAI's ChatGPT-5.5 Pro, which proposed the limit shape in Theorem~\ref{thm:shape} as a conjecture and suggested the coordinatewise transformation $\Phi^{-1}$ as a way to simplify the problem. It suggested the tiling used in Section~\ref{sec:tiling} (for $d>3$) when prompted to extend our tiling idea for $d=3$ to higher dimensions. These suggestions provided starting points for our subsequent development of the argument.

During the development and revision of the paper, we used ChatGPT-5.6 and 6 Pro to help refine proofs, check computations and constants, and identify potential gaps. In particular, AI suggested constants and bounds in the estimates of Lemmas~\ref{lem:constant}, \ref{lem:matching}, and~\ref{lem:cube-lower} to finish our arguments. We evaluated these suggestions, selected which to retain, and revised the resulting arguments and exposition through discussion and review.

We also used AI to help reorganize the manuscript, make systematic changes to notations, and develop Python code for the figures based on our existing visualizations. AI generated an accompanying Lean formalization of the results, available in the repository~\cite{antichain-lean}, for the convenience of some readers. We reviewed the formal Lean statements of the main theorems for consistency with their formulations in this paper, although we have not exhaustively checked all the definitions and results on which they depend. The authors take full responsibility for all results, proofs, references, and exposition.
\printbibliography
\end{document}

%% file: preamble.tex
\usepackage[letterpaper,margin=1in]{geometry}
\usepackage[T1]{fontenc}
\usepackage[utf8]{inputenc}
\usepackage{lmodern,microtype}
\usepackage{amsmath,amssymb,mathtools,amsthm}
\usepackage{bm}
\usepackage[normalem]{ulem}
\usepackage{xcolor,graphicx,enumitem,dsfont,tikz}
\usetikzlibrary{arrows.meta,calc}

\usepackage[titles]{tocloft}

\usepackage{hyperref}
\hypersetup{
  colorlinks=true,
  linkcolor=blue!60!black,
  citecolor=green!40!black,
  urlcolor=blue!70!black,
  linktoc=all,
  bookmarksdepth=2,
  pdftitle={The limit shape of uniform random antichains in [0,1] to the power d},
  pdfsubject={Limit shapes, permutation limits, and applications to coordinate orders}
}

\usepackage[nameinlink,capitalize,noabbrev]{cleveref}
\crefname{section}{Section}{Sections}
\Crefname{section}{Section}{Sections}
\crefname{subsection}{Section}{Sections}
\Crefname{subsection}{Section}{Sections}
\crefname{subsubsection}{Section}{Sections}
\Crefname{subsubsection}{Section}{Sections}

\usepackage[
  backend=biber,
  style=alphabetic,
  maxalphanames=3,
  maxnames=8,
  maxbibnames=20,
  giveninits=true
]{biblatex}

\AtBeginBibliography{\small\sloppy}

\numberwithin{equation}{section}
\newtheorem{theorem}{Theorem}[section]

\newtheorem{proposition}[theorem]{Proposition}
\newtheorem{lemma}[theorem]{Lemma}
\newtheorem{corollary}[theorem]{Corollary}

\theoremstyle{definition}
\newtheorem{question}[theorem]{Question}

\newtheorem{conjecture}[theorem]{Conjecture}

\theoremstyle{remark}
\newtheorem{remark}[theorem]{Remark}

\crefname{theorem}{Theorem}{Theorems}
\crefname{proposition}{Proposition}{Propositions}
\crefname{lemma}{Lemma}{Lemmas}
\crefname{corollary}{Corollary}{Corollaries}
\crefname{conjecture}{Conjecture}{Conjectures}
\crefname{problem}{Problem}{Problems}
\crefname{question}{Question}{Questions}
\Crefname{question}{Question}{Questions}

\setlist{topsep=3pt,itemsep=2pt,parsep=0pt}

%% file: macros.tex
\newcommand{\R}{\mathbb R}

\newcommand{\N}{\mathbb N}

\newcommand{\eps}{\varepsilon}
\newcommand{\dd}{\,\mathrm d}

\renewcommand{\P}{\mathbb P}
\newcommand{\E}{\mathbb E}
\newcommand{\law}{\operatorname{Law}}

\newcommand{\cA}{\mathcal A}
\newcommand{\cB}{\mathcal B}
\newcommand{\cC}{\mathcal C}

\newcommand{\cH}{\mathcal H}

\newcommand{\cS}{\mathcal S}
\newcommand{\cT}{\mathcal T}
\newcommand{\cU}{\mathcal U}
\newcommand{\cV}{\mathcal V}

\newcommand{\cZ}{\mathcal Z}

\newcommand{\ff}{\bm{\mathfrak f}}

\newcommand{\fh}{\mathfrak h}

\newif\iftdo
  \tdotrue

\newcommand{\ind}{\mathds 1}
\newcommand{\vx}{\mathbf x}
\newcommand{\vy}{\mathbf y}

\newcommand{\vf}{\mathbf f}
\newcommand{\vg}{\mathbf g}

\newcommand{\ve}{\mathbf e}
\newcommand{\vv}{\mathbf v}
\newcommand{\vu}{\mathbf u}
\newcommand{\vw}{\mathbf w}
\newcommand{\vX}{\mathbf X}
\newcommand{\vY}{\mathbf Y}
\newcommand{\vG}{\mathbf G}
\newcommand{\vF}{\mathbf F}
\newcommand{\vR}{\mathbf R}

\newcommand{\Ldir}[1]{\mathcal L\left\{#1\right\}}

\DeclareMathOperator{\Vol}{vol}
\DeclareMathOperator{\dist}{dist}
\DeclareMathOperator{\diam}{diam}

\DeclareMathOperator{\Proj}{Proj}
\DeclareMathOperator{\Inv}{Inv}
\DeclareMathOperator{\inv}{inv}